# ON RANDOM ALMOST PERIODIC TRIGONOMETRIC POLYNOMIALS AND APPLICATIONS TO ERGODIC THEORY

By Guy Cohen[1] and Christophe Cuny[2]

*Ben-Gurion University*

We study random exponential sums of the form $\sum_{k=1}^{n} X_k \times \exp\{i(\lambda_k^{(1)} t_1 + \cdots + \lambda_k^{(s)} t_s)\}$, where $\{X_n\}$ is a sequence of random variables and $\{\lambda_n^{(i)} : 1 \leq i \leq s\}$ are sequences of real numbers. We obtain uniform estimates (on compact sets) of such sums, for independent centered $\{X_n\}$ or bounded $\{X_n\}$ satisfying some mixing conditions. These results generalize recent results of Weber [*Math. Inequal. Appl.* **3** (2000) 443–457] and Fan and Schneider [*Ann. Inst. H. Poincaré Probab. Statist.* **39** (2003) 193–216] in several directions. As applications we derive conditions for uniform convergence of these sums on compact sets. We also obtain random ergodic theorems for finitely many commuting measure-preserving point transformations of a probability space. Finally, we show how some of our results allow to derive the Wiener–Wintner property (introduced by Assani [*Ergodic Theory Dynam. Systems* **23** (2003) 1637–1654]) for certain functions on certain dynamical systems.

**1. Introduction.** In their pioneering work, Paley and Zygmund [29] studied Fourier series whose terms have random signs, that is, random Fourier series of the form $\sum_{n=1}^{\infty} \varepsilon_n a_n e^{int}$, where $\{\varepsilon_n\}$ is a Rademacher sequence (i.i.d. random variables taking the values $\pm 1$ with probability $\frac{1}{2}$), and $\{a_n\}$ is a complex sequence. This research was continued by Salem and Zygmund in [34].

Received August 2004; revised February 2005.

[1]Supported in part by Israel Science Foundation Grant 235/01 and by a Skirball Postdoctoral Fellowship of the Center of Advanced Studies in Mathematics of Ben-Gurion University. Also supported by FWF Project P16004–N05 from the E. Schröedinger Institute, Viena.

[2]Supported in part by a Skirball Postdoctoral Fellowship of the Center of Advanced Studies in Mathematics of Ben-Gurion University.

In this paper we obtain uniform estimates of multidimensional random exponential sums of the form $\sum_{k=1}^n X_k e^{i(\lambda_k^{(1)}t_1+\cdots+\lambda_k^{(s)}t_s)}$, where $\{X_n\}$ is a sequence of random variables and $\{\lambda_n^{(i)}:1\leq i\leq s\}$ are sequences of real numbers. Estimations of this kind were obtained (for the one-dimensional case) in [29] and [34], and were extended recently in several directions in [16] and [37].

Such estimates are useful, for instance, to study almost sure (a.s.) uniform convergence of certain random Fourier or random almost periodic series (see, e.g., [12, 16, 20, 34, 37]) and have applications in (random) ergodic theory (e.g., [8]).

Let us be more precise, concerning this latter point.

In the last decades, many authors worked on ergodic theorems with random modulation (sometimes called "randomly weighted ergodic theorems"). One important matter may be formulated as follows: Given a sequence $\{X_n\}$ on a probability space $(\Omega,\mu,\mathcal{F})$, find a measurable set $\Omega^*\subset\Omega$ with $\mu(\Omega^*)=1$, such that for any $\omega\in\Omega^*$ the sequence $a_n:=X_n(\omega)$ is a universally good weight sequence for the ergodic theorem for all functions in some specified class. More precisely, one wants that for any measure-preserving transformation $\tau$ on a probability space $(Y,\Sigma,\pi)$, and any function $f$ on $Y$ with a certain integrability property (e.g., $f\in L_p$), the sequence $\frac{1}{n}\sum_1^n a_k f\circ\tau^k$ converges $\pi$-a.e.

One main tool in the study of such questions (and related ones) is the use of the spectral theorem, which transfers the problem to the study of uniform estimates of random trigonometric polynomials. It seems that the first use of the spectral theorem in this random context appeared in [33], on the base of the results of [34], mentioned above. Then, many authors investigated this direction (see, e.g., [1, 2, 4, 8, 11, 35, 37]).

Another tool, mainly introduced by Rosenblatt [33] in this context (see also the later papers [4] or [7]), is Stein's interpolation theorem, which needs estimates on partial sums of Dirichlet series.

Actually, it seems that what is really needed in order to use Stein's interpolation theorem is to estimate general exponential sums involving Fourier and Dirichlet terms.

This paper may be divided into two parts, estimates and convergence results. First we obtain new estimates, uniform on compacta, for random almost periodic polynomials. Our main result in this direction is the following (see Section 3):

THEOREM 1.1. *Let $\{X_n\}$ be a sequence of random variables, defined on a probability space $(\Omega,\mu)$. Let $\{\lambda_n^{(1)}\},\ldots,\{\lambda_n^{(s)}\}$ be sequences of real numbers.*

(i) *If $\{X_n\}$ are complex-valued, symmetric and independent, then there exist some constants $C,\varepsilon>0$, independent of $\{X_n\}$, such that (with $0/0$*

*interpreted as* 1*)*

$$\left\|\sup_{m\geq 2}\max_{1\leq n<m}\sup_{T\geq 1}\exp\left\{\varepsilon\cdot\left(\max_{(t_1,\ldots,t_s)\in[-T,T]^s}\left|\sum_{k=n+1}^{m}X_k e^{i(\lambda_k^{(1)}t_1+\cdots+\lambda_k^{(s)}t_s)}\right|^2\right)\right.\right.$$

(1) $$\times\left(R_{n,m}\log(T+1)\right.$$

$$\left.\left.\times\log\left(m\vee\max_{1\leq k\leq m}\max_{1\leq i\leq s}|\lambda_k^{(i)}|\right)\right)^{-1}\right\}\right\|_{L_1(\mu)}\leq C,$$

*where* $R_{n,m}=\sum_{k=n+1}^{m}|X_k|^2$. *In particular, for a.e.* $\omega\in\Omega$ *we have*

(2) $$\sup_{m>n}\sup_{T\geq 1}\frac{\max_{(t_1,\ldots,t_s)\in[-T,T]^s}|\sum_{k=n+1}^{m}X_k(\omega)e^{i(\lambda_k^{(1)}t_1+\cdots+\lambda_k^{(s)}t_s)}|^2}{R_{n,m}\log(T+1)\log(m\vee\max_{1\leq k\leq m}\max_{1\leq i\leq s}|\lambda_k^{(i)}|)}<\infty.$$

(ii) *If* $\{X_n\}\subset L_2(\Omega,\mu)$ *are centered independent, then* (1) *and* (2) *remain true with* $R_{n,m}=\sum_{k=n+1}^{m}|X_k|^2+\mathbf{E}|X_k|^2$.

(iii) *If* $\{X_n\}$ *is a bounded martingale difference sequence, then* (1) *and* (2) *remain true with* $R_{n,m}=\sum_{k=n+1}^{m}\|X_k\|_\infty^2$.

We also obtain similar results as in Theorem 1.1 for centered complex bounded random variables which are not necessarily independent. In this case, the quantities $R_{n,m}$ involve some (uniform) correlation coefficients.

This theorem generalizes recent results of Weber [37] and Fan and Schneider [16], the first of which is a one-dimensional version of Theorem 1.1 for periodic polynomials with independent symmetric coefficients, that is, $s=1$ and $\{\lambda_n\}$ is a nondecreasing sequence of natural numbers. The paper of Fan and Schneider gives similar estimations (in a one-dimensional setting) only with $L^1(\mu)$-integrability, while we obtain here an Orlicz space integrability, defined by the function $e^{x^2}-1$. Moreover, Theorem 1.1 shows that Theorem 1 of [16] holds without their (quite restrictive) condition $(\mathcal{V})$. We would like also to underline the fact that in (1) we take the supremum over $T$ inside the integral. This seems to be crucial in the applications such as Stein's interpolation, see Theorem 5.3 below, or for random ergodic theorems for flows that we will explore in a forthcoming paper [6].

Theorem 1.1 is also a generalization of a well-known theorem of Salem and Zygmund [34] for Fourier series whose terms have random signs. Actually, our proof relies on ideas of [34]. It turns out that estimates like (1) are really of a probabilistic nature. In particular, the power of Bernstein's inequality, which was used in [34] (to deal with random trigonometric polynomials), is not needed. It is this remark that allows one to consider general sequences $\{\lambda_n^{(1)}\},\ldots,\{\lambda_n^{(s)}\}$ which are not necessarily integer-valued, and which are not required to satisfy any further assumptions.

To reach the case of general random variables, that is, random variables which are not necessarily independent, we use recent results of [13] (see Proposition 2.5) which permit to deduce exponential inequalities of Azuma's type.

In the second part, we first use our estimates to obtain a.s. uniform convergence on compacta of certain random series of functions of the form $\sum_{n=1}^{\infty} X_n e^{i(\lambda_n^{(1)} t_1 + \cdots + \lambda_n^{(s)} t_s)}$, when $\{X_n\}$ are general bounded or independent (either centered or symmetric). For these results see Section 4.1.

Then we apply our results to ergodic theory. A special case of Theorem 1.1, due to [37], was used in [8] to prove the following:

Let $\{X_n\} \subset L_2(\Omega, \mu)$ be a sequence of centered independent random variables, such that
$$\sum_{n=1}^{\infty} \|X_n\|_2^2 (\log n)^3 < \infty.$$
Then almost surely the sequence $a_n = X_n(\omega)$ has the following property: for any contraction $T$ on $L_2(Y, \pi)$ of a probability space and $f \in L_2(Y, \pi)$, the series $\sum_{n=1}^{\infty} X_n(\omega) T^n f$ converges $\pi$-a.e.

This result raises some questions. If $T$ is induced by a probability-preserving transformation $\tau$, what can be said about functions not in $L_2$, but in some $L_q$, $1 \leq q < 2$? What if we take a sequence of powers $\{j_n\}$? Are there analogues for several commuting transformations? These questions are answered by our main application of Theorem 1.1 in the following (see Section 5):

THEOREM 1.2. *Let $\{X_n\} \subset L_2(\Omega, \mu)$ be a sequence of centered independent random variables, and let $\{j_n^{(1)}\}, \ldots, \{j_n^{(s)}\}$ be sequences of natural numbers. Assume that the series*
$$\sum_{n=1}^{\infty} \|X_n\|_2^2 \log\left(n \vee \max_{1 \leq k \leq n} \max_{1 \leq i \leq s} j_k^{(i)}\right) (\log n)^2$$
*converges. Then there exists a set $\Omega^* \subset \Omega$ of full measure, such that for every $\omega \in \Omega^*$ we have the following: for every commuting family of measure-preserving transformations $\tau_1, \ldots, \tau_s$ on a probability space $(Y, \pi)$, and any $f \in L_q(Y, \pi)$, $1 \leq q \leq 2$, the series*
$$\sum_{n=1}^{\infty} \frac{X_n(\omega) f \circ \tau_1^{j_n^{(1)}} \circ \cdots \circ \tau_s^{j_n^{(s)}}}{n^{(2-q)/2q}}$$
*converges $\pi$-a.e.*

The special case of Theorem 1.2, with $s = 1$, $j_n = n$ and $q = 2$, is also a special case of the result of [8] quoted above (for extensions to the case of two commuting contractions see Theorem 5.2 below).

We conclude the paper by showing connections of our results with the Wiener–Wintner property introduced by Assani in [2]; see Section 6.

**2. Preliminary estimates.**

2.1. *Estimates for almost periodic polynomials.* In [34], Lemma 4.2.3, Salem and Zygmund used Bernstein's inequality to compare the maximum of a trigonometric polynomial with its values on a certain interval. Then they obtained a very sharp result; it seems that in their application, the full strength of their result is not used. We give here some elementary estimates which will be useful in the study of general exponential sums.

LEMMA 2.1. *Let $\{\alpha_n\}$ be a sequence of complex numbers, and let $\{\lambda_n\}$ be a sequence of real numbers. For any $n \geq 1$ put $P_n(t) = \sum_{k=1}^{n} \alpha_k e^{i\lambda_k t}$. Then, for every integer $m \geq 1$ and $T > 0$, we have:*

(i) $\max_{1 \leq n \leq m} \max_{t \in [-T,T]} |P_n'(t)| \leq 3m \max_{1 \leq n \leq m} \{|\lambda_n|\} \cdot \max_{1 \leq n \leq m} \max_{t \in [-T,T]} |P_n(t)|.$

*If in addition $\{\lambda_n\}$ is positive and nondecreasing, then*

(ii) $\max_{1 \leq n \leq m} \max_{t \in [-T,T]} |P_n'(t)| \leq 2\lambda_m \max_{1 \leq n \leq m} \max_{t \in [-T,T]} |P_n(t)|.$

PROOF. Define $P_0 \equiv 0$. For $t \in [-T, T]$ and $1 \leq n \leq m$, we have

$$P_n'(t) = i \sum_{k=1}^{n} \lambda_k \alpha_k e^{i\lambda_k t} = i \sum_{k=1}^{n} \lambda_k (P_k(t) - P_{k-1}(t))$$
$$= i \sum_{k=1}^{n-1} P_k(t)(\lambda_k - \lambda_{k+1}) + i\lambda_n P_n(t).$$

Hence, the results clearly follow. □

REMARKS. 1. When $\{\lambda_n\}$ is a sequence of natural numbers and $T = \pi$, Bernstein's inequality yields that for any $n \geq 1$,

$$\max_{t \in [-\pi,\pi]} |P_n'(t)| \leq \max_{1 \leq k \leq n} \{\lambda_k\} \cdot \max_{t \in [-\pi,\pi]} |P_n(t)|.$$

It is clearly much stronger than (i), and in the monotonic case it is also stronger than (ii).

2. The idea behind the previous simple lemma is to bypass the use of Bernstein's approximation for the derivative of a trigonometric polynomial, in order to overcome the difficulties which appear in trying to extend such an approximation for almost periodic polynomial. Since in this paper the sequence $\{\alpha_n\}$ will always represent a realization of some random variables, it turns out that the obtained result is sufficient for our needs.

NOTATION. For a positive sequence $\{c_n\}$, we define $c_m^* = \max_{1 \leq n \leq m} c_n$.

LEMMA 2.2. *Let $\{\alpha_n\}$ be a sequence of complex numbers, and let $\{\lambda_n^{(1)}\},\ldots,\{\lambda_n^{(s)}\}$ be sequences of real numbers. Put $P_n(t_1,\ldots,t_s) = \sum_{k=1}^n \alpha_k \times e^{i(\lambda_k^{(1)} t_1 + \cdots + \lambda_k^{(s)} t_s)}$. Then for every $T > 0$ and for every integer $m \geq 1$, there exists a rectangle $I \subset [-T,T]^s$, with area $|I|$, satisfying $|I| \geq \prod_{i=1}^s \min\{\frac{1}{6s \cdot m|\lambda_m^{(i)}|^*}, T\}$, such that for every $(t_1,\ldots,t_s) \in I$ we have*

$$\max_{1 \leq n \leq m} |P_n(t_1,\ldots,t_s)| \geq \tfrac{1}{2} \max_{1 \leq n \leq m} \max_{(u_1,\ldots,u_s) \in [-T,T]^s} |P_n(u_1,\ldots,u_s)|.$$

*If in addition $\{\lambda_n^{(1)}\},\ldots,\{\lambda_n^{(s)}\}$ are all positive and nondecreasing, then we can take $I$ with $|I| \geq \prod_{i=1}^s \min\{\frac{1}{4s \cdot \lambda_m^{(i)}}, T\}$.*

PROOF. We first prove the nonmonotonic case. If for some $1 \leq i \leq s$ we have $|\lambda_m^{(i)}|^* = 0$, then the polynomials $\{P_n\}_{n=1}^m$ are all constant with respect to the argument $t_i$. So, we assume that $|\lambda_m^{(i)}|^* \neq 0$ for any $1 \leq i \leq s$.

There exist $1 \leq n_0 \leq m$ and $(u_1,\ldots,u_s) \in [-T,T]^s$, such that

$$M := |P_{n_0}(u_1,\ldots,u_s)| = \max_{1 \leq n \leq m} \max_{(t_1,\ldots,t_s) \in [-T,T]^s} |P_n(t_1,\ldots,t_s)|.$$

Let $(t_1,\ldots,t_s) \in [-T,T]^s$. We have

$$P_{n_0}(t_1,\ldots,t_s) - P_{n_0}(u_1,\ldots,u_s) = \sum_{i=1}^s (t_i - u_i) \frac{\partial P_{n_0}}{\partial t_i}(t_1',\ldots,t_s'),$$

where $(t_1',\ldots,t_s')$ is on the line segment joining $(u_1,\ldots,u_s)$ and $(t_1,\ldots,t_s)$. Hence using Lemma 2.1(i)

$$\begin{aligned}
M &- |P_{n_0}(t_1,\ldots,t_s)| \\
&= ||P_{n_0}(t_1,\ldots,t_s)| - |P_{n_0}(u_1,\ldots,u_s)|| \\
&\leq |P_{n_0}(t_1,\ldots,t_s) - P_{n_0}(u_1,\ldots,u_s)| \\
&\leq \sum_{i=1}^s |t_i - u_i| \left|\frac{\partial P_{n_0}}{\partial t_i}(t_1',\ldots,t_s')\right| \\
&\leq 3m|\lambda_m^{(1)}|^* \max_{1 \leq n \leq m} \max_{v_1 \in [-T,T]} |P_n(v_1, t_2',\ldots,t_s')| \cdot |t_1 - u_1| + \cdots \\
&\quad + 3m|\lambda_m^{(s)}|^* \cdot \max_{1 \leq n \leq m} \max_{v_s \in [-T,T]} |P_n(t_1', t_2',\ldots,t_{s-1}',v_s)| \cdot |t_s - u_s| \\
&\leq 3mM \sum_{i=1}^s |\lambda_m^{(i)}|^* |t_i - u_i|.
\end{aligned}$$

Put

$$I := \left\{(t_1,\ldots,t_s) \in [-T,T]^s : |t_i - u_i| \leq \min\left\{\frac{1}{6s \cdot m|\lambda_m^{(i)}|^*}, T\right\}, i = 1,\ldots,s\right\}.$$

For any $(t_1, \ldots, t_s) \in I$ we have

$$\max_{1 \leq n \leq m} |P_n(t_1, \ldots, t_s)| \geq |P_{n_0}(t_1, \ldots, t_s)| \geq \frac{M}{2}.$$

The monotonic case follows in a similar way. Using Lemma 2.1(ii) we obtain

$$||P_{n_0}(t_1, \ldots, t_s)| - |P_{n_0}(u_1, \ldots, u_s)||$$
$$\leq 2\lambda_m^{(1)} \max_{1 \leq n \leq m} \max_{v_1 \in [-T,T]} |P_n(v_1, t'_2, \ldots, t'_s)| \cdot |t_1 - u_1| + \cdots$$
$$+ 2\lambda_m^{(s)} \max_{1 \leq n \leq m} \max_{v_s \in [-T,T]} |P_n(t'_1, t'_2, \ldots, t'_{s-1}, v_s)| \cdot |t_s - u_s|.$$

In that case we put

$$I := \left\{ (t_1, \ldots, t_s) \in [-T, T]^s : \min\left\{|t_i - u_i| \leq \frac{1}{4s\lambda_m^{(i)}}, T\right\}, i = 1, \ldots, s \right\}. \quad \square$$

DEFINITION 2.1. Let $K$ be a separable compact topological space, and let $\nu$ be a Borel measure on $K$. Let $\{f_n\}$ be a sequence of complex continuous functions on $K$, such that $\sup_{n \geq 1} \max_{x \in K} |f_n(x)| < \infty$. Let $\{\sigma_n\}$ be a nondecreasing sequence, with $\sigma_n \geq 1$. We say that the sequence $\{f_n\}$ forms a $\{\sigma_n\}$-system on $(K, \nu)$ if there exist some constants $\rho_1 > 0$ and $0 < \rho_2 < 1$, such that for every sequence $\{a_n\}$ of complex numbers, and for every $m \geq 1$, there exists a measurable set $I_m \subset K$, with $\nu(I_m) \geq \frac{\rho_1}{\sigma_m}$ such that

$$\max_{1 \leq n \leq m} \left|\sum_{k=1}^n a_k f_k(x)\right| \geq \rho_2 \max_{1 \leq n \leq m} \max_{y \in K} \left|\sum_{k=1}^n a_k f_k(y)\right| \quad \text{for every } x \in I_m.$$

REMARKS. 1. By taking $\{a_n\}$ with zero terms, we can consider the inequality above on any blocks.

2. Definition 2.1 is inspired by the general observation made in [20], Theorem 1, page 68.

EXAMPLE 2.1. Let $\nu$ be the Lebesgue measure on $\mathbf{R}^s$, and let $T \geq 1$. By Lemma 2.2, the sequence $\{e^{i(\lambda_n^{(1)} t_1 + \cdots + \lambda_n^{(s)} t_s)}\}$ forms a $\{(6sn)^s \prod_{i=1}^s (|\lambda_n^{(i)}|^* + 1)\}$-system on $([-T, T]^s, \nu)$. If $\{\lambda_n^{(1)}\}, \ldots, \{\lambda_n^{(s)}\}$ are all positive and nondecreasing, it forms a $\{(4s)^s \prod_{i=1}^s (\lambda_n^{(i)} + 1)\}$-system.

DEFINITION 2.2. Let $K$ be a separable topological space, and let $\nu$ be a Borel measure on $K$. Let $\{K_r\}_{r=1}^\infty \subset K$ be a sequence of compact subspaces of $K$. We say that $\{f_n\}$ forms a *uniform* $\{\sigma_n\}$-system on $\{(K_r, \nu)\}_{r=1}^\infty$ if for every $r \geq 1$ the sequence $\{f_n\}$ forms a $\{\sigma_n\}$-system on $(K_r, \nu)$ with the same corresponding constants $\rho_1$ and $\rho_2$.

EXAMPLE 2.2. For every $r \geq 1$ put $K_r = [-r,r]^s$. Using 2.2 (see also Example 2.1 above), the sequence $\{e^{i(\lambda_n^{(1)} t_1 + \cdots + \lambda_n^{(s)} t_s)}\}$ forms a uniform $\{n^s \prod_{i=1}^s (|\lambda_n^{(i)}|^* + 1)\}$-system on $\{(K_r, \nu)\}_{r=1}^\infty$. If $\{\lambda_n^{(1)}\}, \ldots, \{\lambda_n^{(s)}\}$ are all positive nondecreasing, it forms a uniform $\{\prod_{i=1}^s (\lambda_n^{(i)} + 1)\}$-system.

Let $\{X_n\}$ be a sequence of complex random variables on $(\Omega, \mu)$. Given a sequence $\{f_k\} \subset C(K)$, for any $0 \leq j < l$, we define the random continuous function $S_{j,l} := \sum_{k=j+1}^l X_k f_k$. By separability of $K$ and continuity in $x \in K$ for each fixed $\omega \in \Omega$, we can compute $\|S_{j,l}\| = \max_{x \in K} |\sum_{k=j+1}^l X_k(\omega) f_k(x)|$ as the supremum over a fixed countable dense subset of $K$, so $\|S_{j,l}\|$ is measurable, and $S_{j,l}$ is a $C(K)$-valued random variable.

LEMMA 2.3. *Let $\{X_n\}$ be a sequence of complex random variables on $(\Omega, \mu)$, and let $\{f_n\}$ be a $\{\sigma_n\}$-system on $(K, \nu)$, with constants $\rho_1$ and $\rho_2$ as above. Then for every positive nondecreasing function $\psi$, we have*

$$\mathbf{E}\left(\max_{n<k\leq m} \max_{x\in K} \psi(|S_{n,k}(x)|)\right) \leq \frac{\sigma_m}{\rho_1} \int_K \mathbf{E}\left(\max_{n<k\leq m} \psi\left(\frac{1}{\rho_2}|S_{n,k}(y)|\right)\right)\nu(dy).$$

PROOF. Since $\{f_n\}$ forms a $\{\sigma_n\}$-system on $(K, \nu)$, for every $m \geq 1$, and every $\omega \in \Omega$, there exists a Borel measurable set $I_m(\omega) \subset K$, with $\nu(I_m(\omega)) \geq \rho_1/\sigma_m$, such that for every $x \in K$ we have

$$\psi\left(\max_{n<k\leq m} \max_{y\in K} |S_{n,k}(y)|\right) \mathbf{1}_{I_m(\omega)}(x) \leq \psi\left(\frac{1}{\rho_2} \max_{n<k\leq m} |S_{n,k}(x)|\right).$$

Integrating this inequality on $K$, and using the monotonicity of $\psi$, we obtain

$$\max_{n<k\leq m} \max_{y\in K} \psi(|S_{n,k}(y)|) \leq \frac{\sigma_m}{\rho_1} \int_K \max_{n<k\leq m} \psi\left(\frac{1}{\rho_2}|S_{n,k}(x)|\right)\nu(dx).$$

Taking the expectation yields the result. □

2.2. *Moment inequalities for the partial sums.* We recall here some results that we will use in Section 3. The following lemma is basically Lemma 3 in [29], part I.

LEMMA 2.4. *Let $Z$ be a nonnegative random variable on $(\Omega, \mu)$, and let $C_1$ and $C_2$ be some positive constants. If $\int Z^{2n} d\mu \leq C_1(C_2 n)^n$ for every $n \geq 1$, then $\int \exp(\delta Z^2) d\mu \leq 1 + \frac{C_1}{1-e\delta C_2}$, for every $\delta < \frac{1}{eC_2}$.*

PROOF. Using the estimation $n! \geq \sqrt{2\pi} n^{n+1/2} e^{-n+1/(12n+1)}$ (see [17], page 52), we have

$$\int \exp(\delta Z^2) d\mu = \int \left(1 + \sum_{n=1}^\infty \frac{(\delta Z^2)^n}{n!}\right) d\mu$$

$$\leq 1 + C_1 \sum_{n=1}^{\infty} \frac{(\delta C_2)^n n^n}{n!} \leq 1 + C_1 \sum_{n=1}^{\infty} \frac{(\delta C_2)^n n^n e^n}{n^n \sqrt{n}}$$

$$\leq 1 + C_1 \sum_{n=1}^{\infty} (e\delta C_2)^n \leq 1 + \frac{C_1}{1 - e\delta C_2}. \qquad \square$$

Let $\{X_n\}$ be a sequence of random variables. For any $k \geq 1$, let $\mathcal{F}_k := \sigma(X_1, \ldots, X_k)$ be the $\sigma$-algebra generated by $\{X_1, \ldots, X_k\}$. The following result was obtained by Dedecker.

PROPOSITION 2.5 ([13], Proposition 1). *Let $2 \leq p < \infty$, and let $\{X_n\} \subset L_p(\Omega, \mu)$ be a sequence of real centered random variables. Then we have*

$$(3) \qquad \left\| \sum_{k=1}^{n} X_k \right\|_p \leq \sqrt{2p} \left( \sum_{i=1}^{n} \|X_i^2\|_{p/2} + \sum_{i=2}^{n} \sum_{k=1}^{i-1} \|X_k \mathbf{E}(X_i | \mathcal{F}_k)\|_{p/2} \right)^{1/2}.$$

REMARKS. 1. If $\{X_n\}$ are bounded centered random variables, and in the above inequality we put the $L_\infty$-norm in the right-hand side, then Dedecker's result can be deduced from Theorem 2.4 in [32], page 42.

2. Let $X \in L_r(\Omega, \mu)$, $1 \leq r \leq \infty$, and let $\mathcal{F} \subset \mathcal{F}'$ be two $\sigma$-algebras of $\Omega$. Since the conditional expectation, with respect to a fixed $\sigma$-algebra, contracts the $L_r$-norm, we have

$$(4) \qquad \|\mathbf{E}(X|\mathcal{F})\|_r = \|\mathbf{E}(\mathbf{E}(X|\mathcal{F}')|\mathcal{F})\|_r \leq \|\mathbf{E}(X|\mathcal{F}')\|_r.$$

As a consequence of (4), we may replace $\mathcal{F}_k$ in Dedecker's result by any $\sigma$-algebra $\mathcal{F}_k'$ with respect to which $\{X_1, \ldots, X_k\}$ are measurable. In particular, by usual complexification, it is easy to obtain Dedecker's result for $\{X_n\}$ with complex values. In that case, a factor 2 should be added in front of the right-hand side of (3).

3. As noticed in [13], this result contains Burkholder's inequality for martingale difference sequence.

4. An inspection of the proof of Proposition 1 in [13] shows that the "centered" assumption is not needed.

The following maximal inequality was obtained by Móricz.

PROPOSITION 2.6 ([26], Theorem 1). *Let $1 < p < \infty$, and let $\{X_n\} \subset L_p(\Omega, \mu)$ be a sequence of complex random variables. Assume there exist nonnegative numbers $\{\alpha_n\}$, and some positive constants $C$ and $q > 1$, such that*

$$\left\| \sum_{k=j+1}^{l} X_k \right\|_p^p \leq C \left( \sum_{k=j+1}^{l} \alpha_k \right)^q \qquad \text{for every } l > j \geq 0.$$

*Then for any $m > n \geq 0$,*

$$\left\| \max_{n < l \leq m} \left| \sum_{k=n+1}^{l} X_k \right| \right\|_p^p \leq C_{p,q} \left( \sum_{k=n+1}^{m} \alpha_k \right)^q,$$

where $C_{p,q} = C(1 - \frac{1}{2^{(q-1)/p}})^{-p}$.

REMARKS. 1. Originally, Theorem 1 of [26] is stated for *real* random variables but it extends easily to the case of complex ones. More generally, it extends to the case of Banach-valued random variables.

2. Fix $m > n \geq 0$, and assume that in the assumptions of Proposition 2.6 we only have

$$\left\| \sum_{k=j+1}^{l} X_k \right\|_p^p \leq C \left( \sum_{k=j+1}^{l} \alpha_k \right)^q \qquad \text{for every } n \leq j < l \leq m.$$

This will allow us to estimate the $L_p$-norm of $\max_{n<l\leq m} |\sum_{k=n+1}^{l} X_k|$. Indeed, we put $X'_k = X_k$ and $\alpha'_k = \alpha_k$ for $n < k \leq m$, and otherwise we put $X'_k = 0$ and $\alpha'_k = 0$. Now we apply Proposition 2.6 to $\{X'_k\}$ and $\{\alpha'_k\}$ (see also [8], Proposition 2.3).

3. For fixed $1 < p < \infty$, the quantity $C_{p,q}$ tends monotonically to infinity when $q \downarrow 1$.

4. For $q = 1$, Proposition 2.6 is no longer true. Under the same assumptions, but with $q = 1$, we have (see Theorem 3 in [26] or Proposition 2.2 in [8])

$$\left\| \max_{n < l \leq m} \left| \sum_{k=n+1}^{l} X_k \right| \right\|_p^p \leq C(2 + \log_2 m)^p \sum_{k=n+1}^{m} \alpha_k, \tag{5}$$

for any $m > n \geq 0$.

The proofs of (5) (as given in [26] or [8]) extend to the case of Banach-valued random variables.

**3. Uniform estimates for random polynomials.** Let $\{X_n\} \subset L_\infty(\Omega, \mu)$ be a sequence of complex centered random variables and let $\{f_k\} \subset C(K)$. For any $m > n \geq 0$ define $S_{n,m}(x) := \sum_{k=n+1}^{m} X_k f_k(x)$ and

$$R_{n,m} := \sum_{i=n+1}^{m} \|X_i\|_\infty^2 + \sum_{i=n+2}^{m} \sum_{k=1}^{i-1} \|X_k \mathbf{E}(X_i | \mathcal{F}_k)\|_\infty. \tag{6}$$

THEOREM 3.1. *Let $\{X_n\} \subset L_\infty(\Omega, \mu)$ be a sequence of complex centered random variables. Let $\{f_n\}$ be a $\{\sigma_n\}$-system on $(K, \nu)$, with corresponding constants $\rho_1$ and $\rho_2$. Then for every $m > n \geq 0$ we have (with $0/0$ interpreted as $1$)*

$$\left\| \max_{n<l\leq m} \exp\left\{ \varepsilon \cdot \frac{\max_{x\in K} |\sum_{k=n+1}^{l} X_k f_k(x)|^2}{R_{n,m}} \right\} \right\|_{L_1(\mu)} \leq \frac{3\nu(K)\sigma_m}{\rho_1},$$

*where $\varepsilon = \rho_2^2/(6400\,\mathbf{e} \cdot \max_{n<k\leq m} \max_{x\in K} |f_k(x)|^2)$. In particular, for every $1 \leq p < \infty$ we have*

$$\left\| \max_{n<l\leq m} \max_{x\in K} \left| \sum_{k=n+1}^{l} X_k f_k(x) \right| \right\|_{L_p(\mu)} \leq C_p \sqrt{R_{n,m} \log(\sigma_m + 1)},$$

*where $C_p = \sqrt{[1 + 2\log(e^{p/2} + 3\nu(K)/\rho_1)]/\varepsilon}$.*

PROOF. By definition, $\{f_n\}$ is uniformly bounded. By homogeneity, it is enough to prove the theorem for the case where $\sup_{n\geq 1} \max_{x\in K} |f_k(x)| \leq 1$. Fix $x \in K$, and let $p \geq 2$. By application of Proposition 2.5 (see remarks 1 and 2 after its formulation) to the sequence $\{X_k f_k(x)\} \subset L_\infty(\Omega, \mu)$ we have $\|S_{n,m}(x)\|_p \leq \sqrt{8pR_{n,m}(x)}$, where

$$R_{n,m}(x) := \sum_{i=n+1}^{m} \|X_i f_i(x)\|_\infty^2 + \sum_{i=n+2}^{m} \sum_{k=1}^{i-1} \|X_k f_k(x) \mathbf{E}(X_i f_i(x) | \mathcal{F}_k(x))\|_\infty,$$

(1)
and $\mathcal{F}_k(x) := \sigma(X_1 f_1(x), \ldots, X_k f_k(x))$. Clearly, we have $\mathcal{F}_k(x) \subset \mathcal{F}_k = \sigma(X_1, \ldots, X_k)$ for any $x \in K$. Since $X_k f_k(x)$ is measurable with respect to $\mathcal{F}_k(x)$, we may and do apply inequality (4) to obtain (1) with $\mathcal{F}_k(x)$ replaced by $\mathcal{F}_k$. Using the assumption $|f_n(x)| \leq 1$, we have $R_{n,m}(x) \leq R_{n,m}$, where $R_{n,m}$ is defined by (6). Hence for every $x \in K$ and $p \geq 2$ we have

$$\|S_{n,m}(x)\|_p \leq \sqrt{8pR_{n,m}},$$

for every $m > n \geq 0$.

For every $i \geq 1$ put

$$\alpha_i = \sum_{k=1}^{i} \|X_k \mathbf{E}(X_i | \mathcal{F}_k)\|_\infty.$$

So, $R_{n,m} = \sum_{i=n+1}^{m} \alpha_i$. Fix $p_0 > 2$, and take any $p \geq p_0 > 2$. By Proposition 2.6, applied with $q = p/2 > 1$ and $\{\alpha_i\}$, we obtain

$$\left\| \max_{n<l\leq m} |S_{n,l}(x)| \right\|_p \leq (1 - 2^{(1-p/2)/p})^{-1} \sqrt{8pR_{n,m}}$$

$$\leq (1 - 2^{(1-p_0/2)/p_0})^{-1} \sqrt{8pR_{n,m}}.$$

For $2 \leq p \leq p_0$, we use $\|\cdot\|_p \leq \|\cdot\|_{p_0}$. Finally, we obtain for every $2 \leq p < \infty$ and $x \in K$

$$\left\|\max_{n < l \leq m} |S_{n,l}(x)|\right\|_p \leq \sqrt{8pC_{p_0} R_{n,m}},$$

where we can take $C_{p_0} = p_0(1 - 2^{(1-p_0/2)/p_0})^{-2}$.

By application of Lemma 2.4, with $C_1 = 1$ and $C_2 = 16C_{p_0} R_{n,m}$, we deduce

$$\int \max_{n < l \leq m} \exp(\delta |S_{n,l}(x)|^2) \, d\mu = \int \exp\left(\delta \max_{n < l \leq m} |S_{n,l}(x)|^2\right) d\mu$$

$$\leq 1 + \frac{1}{1 - 16e\delta C_{p_0} R_{n,m}},$$

for every $\delta < \frac{1}{16eC_{p_0} R_{n,m}}$, and for every $x \in K$.

Let $\rho_1$ and $\rho_2$ be the corresponding constants for the $\{\sigma_n\}$-system $\{f_n\}$, as given in Definition 2.1. By Lemma 2.3, applied to the function $u \mapsto \exp(\delta u^2)$ for some $\delta < \frac{\rho_2^2}{16eC_{p_0} R_{n,m}}$, we deduce

$$\int_\Omega \max_{n < l \leq m} \exp\left(\delta \max_{x \in K} |S_{n,l}(x)|^2\right) d\mu$$

$$\leq \frac{\sigma_m}{\rho_1} \int_K \int_\Omega \max_{n < l \leq m} \exp\left(\frac{\delta}{\rho_2^2} |S_{n,l}(x)|^2\right) d\mu \, d\nu(x)$$

$$\leq \frac{\nu(K)\sigma_m}{\rho_1}\left(1 + \frac{1}{1 - 16e(\delta/\rho_2^2)C_{p_0} R_{n,m}}\right).$$

Put $\delta = \rho_2^2/(32\mathbf{e} C_{p_0} R_{n,m})$ and $\varepsilon_{p_0} = \rho_2^2/(32\mathbf{e} C_{p_0})$. Since $p_0 > 2$ is arbitrary, and $\min_{p_0 \geq 2}\{C_{p_0}\} < 200$, we can choose $p_0 > 2$ such that $\varepsilon = \varepsilon_{p_0} = \rho_2^2/(6400\mathbf{e})$. This yields the first result.

For $p > 0$, define the function $\phi_p(x) = (\log(e^p + x))^p$ for any $x \geq 0$. Then $\phi_p$ is concave, and $(\log x)^p \leq \phi_p(x)$ for $x \geq 1$. Hence, by Jensen's inequality and the first result, we obtain

$$\mathbf{E}\left\{\varepsilon^{p/2} \cdot \frac{(\max_{n < l \leq m} \max_{x \in K} |S_{n,l}|)^p}{R_{n,m}^{p/2}}\right\}$$

$$\leq \mathbf{E}\left\{\phi_{p/2}\left(\exp\left\{\varepsilon \cdot \frac{\max_{n < l \leq m} \max_{x \in K} |S_{n,l}|^2}{R_{n,m}}\right\}\right)\right\}$$

$$\leq \phi_{p/2}\left(\frac{3\nu(K)\sigma_m}{\rho_1}\right) = \left(\log\left(e^{p/2} + \frac{3\nu(K)\sigma_m}{\rho_1}\right)\right)^{p/2}.$$

Hence,

$$\left\|\max_{n<l\leq m}\max_{x\in K}|S_{n,l}|\right\|_p^2 \leq \frac{R_{n,m}}{\varepsilon}\log((e^{p/2}+3\nu(K)/\rho_1)\sigma_m)$$

$$\leq \frac{R_{n,m}}{\varepsilon}\left(1+\frac{\log(e^{p/2}+3\nu(K)/\rho_1)}{\log(\sigma_m+1)}\right)\log(\sigma_m+1).$$

Since $\sigma_m \geq 1$, we obtain the second result with $C_p$ as defined in the statement. $\square$

REMARKS. 1. The above theorem generalizes ideas of Theorem 4.3.1 and Lemma 5.1.3 in [34] (see also [38], Chapter V, Theorem 8.34 and [20], Theorem 1, page 69).

2. Clearly, the first statement of the theorem yields the finiteness of all the $L_p$-norms. In fact, as noticed in [22], Lemma 3.7, page 66, it yields that all the $L_p$-norms, $2 \leq p < \infty$, are comparable to the $L_2$-norm. This in turn implies the comparability of all $L_p$-norms, with $p \geq 2$. This is a weak generalization of the Kahane–Khinchine inequality for Rademacher series (see [23], Corollary 4.6, page 43 and see also [20], page 282).

3. It is clear from Definition 2.1 that constant functions, that is, $f_n(x) \equiv C \neq 0$ for every $n \geq 1$, form a $\{\sigma_n\}$-system for any nondecreasing $\{\sigma_n\}$. By doing this, we may obtain classical maximal inequalities for sums of bounded random variables.

COROLLARY 3.2. *Let $1 < p \leq 2$, and let $\{X_n\} \subset L_p(\Omega,\mu)$ be a sequence of complex centered independent random variables. Let $\{f_n\}$ be a $\{\sigma_n\}$-system on $(K,\nu)$, with corresponding constants $\rho_1$ and $\rho_2$. Then for every $m > n \geq 0$ we have*

$$\left\|\max_{n<l\leq m}\max_{x\in K}\left|\sum_{k=n+1}^{l}X_kf_k(x)\right|\right\|_{L_p(\mu)} \leq 2C_p\sqrt{\log(\sigma_m+1)}\left(\sum_{k=n+1}^{m}\|X_k\|_p^p\right)^{1/p},$$

*where $C_p = \sqrt{[1+2\log(e^{p/2}+3\nu(K)/\rho_1)]/\varepsilon}$.*

PROOF. We first prove the corollary when $\{X_n\}$ are symmetric random variables. Let $\{\varepsilon_n\}$ be a Rademacher sequence, which is independent of $\{X_n\}$, and let $\mathbf{E}_\varepsilon$ be the corresponding expectation with respect to the probability space of $\{\varepsilon_n\}$. Let $\mathbf{E}$ be the expectation in $(\Omega,\mu)$. Fix $\omega \in \Omega$. We apply the second result of Theorem 3.1 to the independent sequence $\{X_n(\omega)\varepsilon_n\}$, in order to obtain

$$\mathbf{E}_\varepsilon\left[\max_{n<l\leq m}\max_{x\in K}\left|\sum_{k=n+1}^{l}X_k(\omega)\varepsilon_kf_k(x)\right|^p\right]$$

$$\leq (C_p/2)^p \left(\log(\sigma_m + 1) \sum_{k=n+1}^{m} |X_k(\omega)\varepsilon_k|^2\right)^{p/2}$$

$$\leq (C_p/2)^p (\log(\sigma_m + 1))^{p/2} \sum_{k=n+1}^{m} |X_k(\omega)|^p,$$

where $C_p$ may be computed from Theorem 3.1.

By taking the expectation $\mathbf{E}$ in the above inequality, we obtain that

$$\mathbf{EE}_\varepsilon \left[\max_{n<l\leq m} \max_{x\in K} \left|\sum_{k=n+1}^{l} X_k \varepsilon_k f_k(x)\right|^p\right]$$

$$\leq (C_p/2)^p (\log(\sigma_m + 1))^{p/2} \sum_{k=n+1}^{m} \mathbf{E}|X_k|^p.$$

Since $\{X_n\}$ is symmetric and independent, and is also independent of $\{\varepsilon_n\}$, the sequences $\{X_n\varepsilon_n\}$ and $\{X_n\}$ are stochastically equivalent. So, we have

$$\mathbf{E}\left[\max_{n<l\leq m} \max_{x\in K} \left|\sum_{k=n+1}^{l} X_k f_k(x)\right|^p\right]$$

$$\leq (C_p/2)^p (\log(\sigma_m + 1))^{p/2} \sum_{k=n+1}^{m} \mathbf{E}|X_k|^p.$$

This proves the corollary for the symmetric case.

Now, let $\{X'_n\}$ be an independent copy of $\{X_n\}$, defined on $(\Omega', \mu')$, and let $\mathbf{E}'$ be the corresponding expectation. Clearly, the sequence $\{X_n - X'_n\}$ is a symmetric sequence, so we apply the first result of the proof, for the symmetric case, in order to obtain

$$\mathbf{EE}'\left[\max_{n<l\leq m} \max_{x\in K} \left|\sum_{k=n+1}^{l} (X_k - X'_k) f_k(x)\right|^p\right]$$

$$\leq (C_p/2)^p (\log(\sigma_m + 1))^{p/2} \sum_{k=n+1}^{m} \mathbf{EE}'|X_k - X'_k|^p.$$

Using Jensen's inequality and the fact that $\mathbf{E}'(X'_k) = 0$, we have

$$\mathbf{EE}'\left[\max_{n<l\leq m} \max_{x\in K} \left|\sum_{k=n+1}^{l} (X_k - X'_k) f_k(x)\right|^p\right]$$

$$\geq \mathbf{E}\left[\max_{n<l\leq m} \max_{x\in K} \mathbf{E}'\left|\sum_{k=n+1}^{l} (X_k - X'_k) f_k(x)\right|^p\right]$$

$$\geq \mathbf{E}\left[\max_{n<l\leq m}\max_{x\in K}\left|\sum_{k=n+1}^{l}X_k f_k(x)\right|^p\right].$$

Hence,

$$\mathbf{E}\left[\max_{n<l\leq m}\max_{x\in K}\left|\sum_{k=n+1}^{l}X_k f_k(x)\right|^p\right]$$

$$\leq (C_p/2)^p (\log(\sigma_m+1))^{p/2} \sum_{k=n+1}^{m} \mathbf{E}'\mathbf{E}'|X_k - X_k'|^p$$

$$\leq 2^p (C_p/2)^p (\log(\sigma_m+1))^{p/2} \sum_{k=n+1}^{m} \mathbf{E}|X_k|^p,$$

and the result follows. □

NOTATION. Let $s \geq 1$, and consider the Euclidean space $\mathbf{R}^s$. We denote by boldface, for example, $\mathbf{t} = (t_1, \ldots, t_s)$, a vector in $\mathbf{R}^s$. For any $\mathbf{t}, \mathbf{u} \in \mathbf{R}^s$ we denote by $\langle \mathbf{t}, \mathbf{u} \rangle = t_1 u_1 + \cdots + t_s u_s$ the inner product in $\mathbf{R}^s$. We recall our notation $c_m^* = \max_{1\leq n\leq m} c_n$, for any positive sequence $\{c_n\}$.

COROLLARY 3.3. *Let $1 < p \leq 2$, and let $\{X_n\} \subset L_p(\Omega, \mu)$ be a sequence of complex centered independent random variables. Let $\boldsymbol{\lambda}_n = (\lambda_n^{(1)}, \ldots, \lambda_n^{(s)})$ be a sequence of vectors in $\mathbf{R}^s$, and let $T \geq 1$. Then there exists a positive constant $C_p$, which does not depend on $\{X_n\}$, such that for every $m > n \geq 0$ we have*

$$\left\|\max_{n<l\leq m}\max_{\mathbf{t}\in[-T,T]^s}\left|\sum_{k=n+1}^{l}X_k e^{i\langle \boldsymbol{\lambda}_n, \mathbf{t}\rangle}\right|\right\|_{L_p(\mu)}$$

$$\leq C_p \sqrt{\log\left[m^s \prod_{i=1}^{s}(|\lambda_m^{(i)}|^* + 1) + 1\right]}\left(\sum_{k=n+1}^{m}\|X_k\|_p^p\right)^{1/p}.$$

PROOF. Let $K = [-T, T]^s$, and let $\nu$ be the Lebesgue measure on $\mathbf{R}^s$. By Example 2.1 the sequence $\{e^{i\langle \boldsymbol{\lambda}_n, \mathbf{t}\rangle}\}$ forms a $\{m^s \prod_{i=1}^{s}(|\lambda_m^{(i)}|^* + 1)\}$-system on $(K, \nu)$. Colorally 3.2 yields the result. □

Now we present further applications of Lemma 2.3. Let $\{X_n\} \subset L_p(\mu)$, $2 \leq p < \infty$, be centered random variables. For any $m > n \geq 0$ define

$$R_{n,m}^{(p)} = \sum_{i=n+1}^{m}\|X_i\|_p^2 + \sum_{i=n+2}^{m}\sum_{k=1}^{i-1}\|X_k \mathbf{E}(X_i|\mathcal{F}_k)\|_{p/2}.$$

In the case of unbounded random variables, we can say the following:

THEOREM 3.4. *Let $2 < p < \infty$, and let $\{X_n\} \subset L_p(\Omega, \mu)$ be a sequence (not necessarily independent) of complex centered random variables. Let $\{f_n\}$ be a $\{\sigma_n\}$-system on $(K, \nu)$, with corresponding constants $\rho_1$ and $\rho_2$. Then for every $m > n \geq 0$ we have*

$$\left\| \max_{n < l \leq m} \max_{x \in K} \left| \sum_{k=n+1}^{l} X_k f_k(x) \right| \right\|_p \leq C_p \max_{n < k \leq m} \max_{x \in K} |f_k(x)| (\sigma_m)^{1/p} \sqrt{R_{n,m}^{(p)}},$$

*where $C_p = \frac{2\sqrt{2p}}{(\rho_1)^{1/p} \rho_2}(1 - 2^{(1-p/2)/p})^{-1}$.*

PROOF. The proof starts as the proof of Theorem 3.1. We use Dedecker's inequality to obtain

$$\|S_{n,m}(x)\|_p \leq \max_{n < k \leq m} \max_{x \in K} |f_k(x)| \cdot \sqrt{8p R_{n,m}^{(p)}},$$

for every $m > n \geq 0$. Then, using Móricz's inequality [26] we obtain

$$\left\| \max_{n < l \leq m} |S_{n,l}(x)| \right\|_p \leq \max_{n < k \leq m} \max_{x \in K} |f_k(x)| (1 - 2^{(1-p/2)/p})^{-1} \sqrt{8p R_{n,m}^{(p)}}.$$

By Lemma 2.3, applied to the function $u \mapsto u^p$, we obtain the result. $\square$

REMARKS. 1. Using remark 4 after Proposition 2.6, for $p = 2$ we obtain

$$\left\| \max_{n < l \leq m} \max_{x \in K} \left| \sum_{k=n+1}^{l} X_k f_k(x) \right| \right\|_2$$
$$\leq C_2 \log(4m) \max_{n < k \leq m} \max_{x \in K} |f_k(x)| (\sigma_m)^{1/2} \sqrt{R_{n,m}^{(2)}},$$

where $C_2 = \frac{4}{\sqrt{\rho_1 \rho_2 \log 2}}$.

2. Let $\{j_k\}$ be a sequence of integers, and let $f_k(t) = e^{ij_k t}$. Using the Cauchy–Schwarz inequality one can see that for every $\omega \in \Omega$ we have

$$\max_{n < l \leq m} \max_{t \in [-\pi, \pi)} \left| \sum_{k=n+1}^{l} X_k(\omega) e^{ij_k t} \right| \leq C \sqrt{m \sum_{k=n+1}^{m} |X_k(\omega)|^2}.$$

This illustrates the limitation of the method when $p = 2$.

THEOREM 3.5. *Let $\{X_n\} \subset L_\infty(\Omega, \mu)$ be a sequence of complex centered random variables. Let $\{f_n\}$ be a uniform $\{\sigma_n\}$-system on $\{(K_r, \nu)\}$. Assume that for some $q \geq 1$, we have $\{1/\nu(K_r)\} \in \ell_{2q}$ and $\{1/\sigma_n\} \in \ell_q$. Then there exists some positive constants $\varepsilon$ and $C$, independent of $\{X_n\}$, such that (with $0/0$ interpreted as $1$)*

$$\left\| \sup_{m \geq 1} \max_{0 \leq n < m} \sup_{r \geq 1} \exp \left\{ \varepsilon \cdot \frac{\max_{x \in K_r} |\sum_{k=n+1}^{m} X_k f_k(x)|^2}{R_{n,m} \log(\nu(K_r)\sigma_m + 1)} \right\} \right\|_{L_1(\mu)} \leq C.$$

*Hence for a.e.* $\omega \in \Omega$ *we have*

$$\sup_{m>n} \sup_{r\geq 1} \frac{\max_{x\in K_r}|\sum_{k=n+1}^m X_k(\omega)f_k(x)|^2}{R_{n,m}\log(\nu(K_r)\sigma_m + 1)} < \infty.$$

PROOF. For any $\omega \in \Omega$ put $S_{n,m}(\omega)(x) = \sum_{k=n+1}^m X_k(\omega)f_k(x)$, and when $\omega$ is not specified put $S_{n,m}(x)$. By assumptions $\{\nu(K_r)\}$ and $\{\sigma_n\}$ tend to infinity, so without loss of generality we may and do assume that $\log(\nu(K_r)\sigma_m + 1) \geq 1$ for every $m, r \geq 1$. This assumption reflects only in the values of $\varepsilon$ and $C$.

The sequence $\{f_n\}$ forms a *uniform* $\{\sigma_n\}$-system, hence for every $r \geq 1$, the sequence $\{f_n\}$ forms a $\{\sigma_n\}$-system on $(K_r, \nu)$, with corresponding constants $\rho_1$ and $\rho_2$. For every $r \geq 1$, we apply Theorem 3.1 with $(K_r, \nu)$ and $\rho_1, \rho_2$. Hence, there exist *universal* constants $C_1, \varepsilon > 0$, such that for every $m > n \geq 0$ we have

$$\left\|\exp\left\{\varepsilon \cdot \frac{\max_{x\in K_r}|S_{n,m}(x)|^2}{R_{n,m}}\right\}\right\|_{L_1(\mu)} \leq C_1 \nu(K_r)\sigma_m,$$

for any $r \geq 1$. Hence for any $m > n \geq 0$ and $r \geq 1$, we have

(1)
$$\left\|\exp\left\{\varepsilon \cdot \frac{\max_{x\in K_r}|S_{n,m}(x)|^2}{R_{n,m}} - (2q+1)\log(\nu(K_r)\sigma_m)\right\}\right\|_{L_1(\mu)}$$
$$\leq C_1(\nu(K_r)\sigma_m)^{-2q}.$$

For any $m > n \geq 0$ and $r \geq 1$ put

$$I_{n,m,r} = \left\{\omega \in \Omega : \varepsilon \cdot \max_{x\in K_r}|S_{n,m}(x)(\omega)|^2 \geq (2q+1)R_{n,m}\log(\nu(K_r)\sigma_m + 1)\right\}.$$

Using (1) we have

$$\left\|\exp\left\{\varepsilon \cdot \frac{\max_{x\in K_r}|S_{n,m}(x)|^2}{R_{n,m}\log(\nu(K_r)\sigma_m + 1)} - (2q+1)\right\}\mathbf{1}_{I_{n,m,r}}\right\|_{L_1(\mu)}$$

$$\leq \left\|\exp\left\{\log(\nu(K_r)\sigma_m + 1)\right.\right.$$
$$\left.\left.\times\left[\varepsilon \cdot \frac{\max_{x\in K_r}|S_{n,m}(x)|^2}{R_{n,m}\log(\nu(K_r)\sigma_m + 1)} - (2q+1)\right]\right\}\mathbf{1}_{I_{n,m,r}}\right\|_{L_1(\mu)}$$

$$\leq \frac{C_1}{(\nu(K_r)\sigma_m)^{2q}}.$$

Hence

$$\sum_{n=0}^{m-1}\sum_{r\geq 1}\left\|\exp\left\{\varepsilon \cdot \frac{\max_{x\in K_r}|S_{n,m}(x)|^2}{R_{n,m}\log(\nu(K_r)\sigma_m + 1)} - (2q+1)\right\}\mathbf{1}_{I_{n,m,r}}\right\|_{L_1(\mu)}$$

$$\leq \frac{mC_1 \|1/\nu(K_r)\|_{\ell_{2q}}^{2q}}{(\sigma_m)^{2q}}.$$

By assumption, the nonincreasing sequence $\{1/\sigma_n\}$ is in $\ell_q$. Hence by Kronecker's lemma, $\sup_{m\geq 1}\{\frac{m}{\sigma_m^q}\} < \infty$. We deduce

$$\sum_{m=1}^{\infty} \sum_{n=0}^{m-1} \sum_{r\geq 1} \left\| \exp\left\{\varepsilon \cdot \frac{\max_{x\in K_r}|S_{n,m}(x)|^2}{R_{n,m}\log(\nu(K_r)\sigma_m+1)} - (2q+1)\right\} \mathbf{1}_{I_{n,m,r}} \right\|_{L_1(\mu)} \leq C_2,$$

where $C_2 = C_1 \|1/\nu(K_r)\|_{\ell_{2q}}^{2q} \cdot \sup_{m\geq 1}\{\frac{m}{\sigma_m^q}\} \cdot \|\{\frac{1}{\sigma_n}\}\|_{\ell_q}^q$. So,

$$\left\| \sup_{m\geq 1} \max_{0\leq n<m} \sup_{r\geq 1} \exp\left\{\varepsilon \cdot \frac{\max_{x\in K_r}|S_{n,m}(x)|^2}{R_{n,m}\log(\nu(K_r)\sigma_m+1)} - (2q+1)\right\} \mathbf{1}_{I_{n,m,r}} \right\|_{L_1(\mu)} \leq C_2.$$

But, if $\omega \notin I_{m,n,r}$ for some $m > n \geq 0$ and $r \geq 1$, then

$$\varepsilon \cdot \frac{\max_{x\in K_r}|S_{n,m}(x)|^2}{R_{n,m}\log(\nu(K_r)\sigma_m+1)} \leq 2q+1.$$

So,

$$(7) \quad \left\| \sup_{m\geq 1} \max_{0\leq n<m} \sup_{r\geq 1} \exp\left\{\varepsilon \cdot \frac{\max_{x\in K_r}|S_{n,m}(x)|^2}{R_{n,m}\log(\nu(K_r)\sigma_m+1)}\right\} \right\|_{L_1(\mu)} \leq C_2 e^{2q+1}.$$

This proves the first assertion, so in particular the integrand is finite a.e., which is the second statement. □

REMARKS. 1. The quantity $\varepsilon$ and $C$ in the above theorem can be completely computed by Theorem 3.1. Note that $C_1 = 3/\rho_1$ in the proof of Theorem 3.5.

2. The conditions $\{1/\sigma_n\} \in \ell_r$, for some $r > 0$, together with nondecreasingness of $\{\sigma_n\}$ (by definition) is equivalent by Kronecker's lemma to $\sigma_n \geq Cn^\delta$, for some $\delta > 0$ and for every $n \geq 1$. Since the assertions of Therem 3.5 are not affected by removing finite numbers of pairs $\{m > n\}$, we could replace the assumption $\{1/\sigma_n\} \in \ell_{2q}$ by $\sigma_n \geq Cn^\delta$, for some $\delta > 0$ and for every $n \geq 1$.

3. If $\{1/\sigma_n\}$ is not in $\ell_r$ for any $r > 0$, then $\{f_n\}$ is still a $\{\max\{n,\sigma_n\}\}$-system. With respect to this system, $C$ in the above theorem is independent of $\{\sigma_n\}$.

COROLLARY 3.6. *Let $\{X_n\}$ be symmetric independent complex valued random variables on $(\Omega,\mu)$. Let $\{f_n\}$ be a uniform $\{\sigma_n\}$-system on $\{(K_r,\nu)\}$. Assume that for some $q \geq 1$, we have $\{1/\nu(K_r)\} \in \ell_{2q}$ and $\{1/\sigma_n\} \in \ell_q$. Then* *there exist some positive constants $\varepsilon$ and $C$, independent of $\{X_n\}$, such that (with $0/0$ interpreted as 1)*

$$\left\| \sup_{m\geq 1} \max_{0\leq n<m} \sup_{r\geq 1} \exp\left\{ \varepsilon \cdot \frac{\max_{x\in K_r} |\sum_{k=n+1}^m X_k f_k(x)|^2}{\log(\nu(K_r)\sigma_m + 1)\sum_{k=n+1}^m |X_k|^2} \right\} \right\|_{L_1(\mu)} \leq C.$$

PROOF. Let $\{\varepsilon_n\}$ be a Rademacher sequence which is independent of $\{X_n\}$. Let $\mathbf{E}_\varepsilon$ and $\mathbf{E}$ be the corresponding expectations in the probability spaces of $\{\varepsilon_n\}$ and $\{X_n\}$, respectively. For a.e. $\omega \in \Omega$, the sequence $\{X_n(\omega)\varepsilon_n\}$ is a sequence of independent bounded random variables. So, we may apply Theorem 3.5. Hence there exist some positive constants $\varepsilon$ and $C$, which are *independent* of $\{X_n(\omega)\varepsilon_n\}$, such that for a.e. $\omega \in \Omega$ we have

$$\mathbf{E}_\varepsilon \left[ \sup_{m\geq 1} \max_{0\leq n<m} \sup_{r\geq 1} \exp\left\{ \varepsilon \cdot \frac{\max_{x\in K_r} |\sum_{k=n+1}^m X_k(\omega)\varepsilon_k f_k(x)|^2}{\log(\nu(K_r)\sigma_m + 1)\sum_{k=n+1}^m |X_k(\omega)\varepsilon_k|^2} \right\} \right] \leq C.$$

By taking the expectation $\mathbf{E}$ we have

$$\mathbf{E}\mathbf{E}_\varepsilon \left[ \sup_{m\geq 1} \max_{0\leq n<m} \sup_{r\geq 1} \exp\left\{ \varepsilon \cdot \frac{\max_{x\in K_r} |\sum_{k=n+1}^m X_k\varepsilon_k f_k(x)|^2}{\log(\nu(K_r)\sigma_m + 1)\sum_{k=n+1}^m |X_k\varepsilon_k|^2} \right\} \right] \leq C.$$

Since $\{X_n\}$ is symmetric and independent, and is also independent of $\{\varepsilon_n\}$, the sequences $\{X_n\varepsilon_n\}$ and $\{X_n\}$ are stochastically equivalent. Hence, the assertion of the theorem follows from the above result for $\{X_n\varepsilon_n\}$. □

For the next results, we will specify the $\{\sigma_n\}$-system.

COROLLARY 3.7. *Let $\{X_n\}$ be symmetric independent complex valued random variables on $(\Omega, \mu)$. Let $\boldsymbol{\lambda}_n = (\lambda_n^{(1)}, \ldots, \lambda_n^{(s)})$ be a sequence of vectors in $\mathbf{R}^s$. Then there exist some positive constants $\varepsilon$ and $C$, independent of $\{X_n\}$, such that (with $0/0$ interpreted as 1)*

$$\left\| \sup_{m\geq 1} \max_{0\leq n<m} \sup_{T\geq 1} \exp\left\{ \varepsilon \cdot \left( \max_{\mathbf{t}\in[-T,T]^s} \left| \sum_{k=n+1}^m X_k e^{i\langle \boldsymbol{\lambda}_k, \mathbf{t}\rangle} \right|^2 \right) \right.\right.$$

$$\times \left( \log\left[ C_s T^s \cdot m^s \cdot \prod_{i=1}^s (|\lambda_m^{(i)}|^* + 1) + 1 \right] \right.$$

$$\left.\left.\left. \times \sum_{k=n+1}^m |X_k|^2 \right)^{-1} \right\} \right\|_{L_1(\mu)} \leq C,$$

*where $C_s = (12s)^s$. Hence for a.e. $\omega \in \Omega$ we have*

$$\sup_{m>n} \sup_{T\geq 1} \frac{\max_{\mathbf{t}\in[-T,T]^s} |\sum_{k=n+1}^m X_k(\omega) e^{i\langle \boldsymbol{\lambda}_k, \mathbf{t}\rangle}|^2}{\log[C_s T^s \cdot m^s \cdot \prod_{i=1}^s (|\lambda_m^{(i)}|^* + 1) + 1] \cdot \sum_{k=n+1}^m |X_k(\omega)|^2} < \infty.$$

PROOF. For every $T \geq 1$ and $\mathbf{t} \in \mathbf{R}^s$, put $K_T = [-T,T]^s$ and $f_n(\mathbf{t}) = e^{i\langle \boldsymbol{\lambda}_k, \mathbf{t}\rangle}$. By uniform continuity, the measurable function $\max_{\mathbf{t} \in [-T,T]^s} \times |\sum_{k=n+1}^m X_k e^{i\langle \boldsymbol{\lambda}_k, \mathbf{t}\rangle}|$ is a continuous function of $T$. Hence, the suprema over $T \geq 1$ can be taken as suprema over the rational numbers. So, the integrand is measurable.

As noted in Example 2.2, the sequence of functions $f_n(\mathbf{t}) = e^{i\langle \boldsymbol{\lambda}_k, \mathbf{t}\rangle}$ forms a uniform $\{(6s)^s m^s \prod_{i=1}^s (|\lambda_n^{(i)}|^* + 1)\}$-system on $\{(K_r, \nu)\}_{r=1}^\infty$, where $\nu$ is the Lebesgue measure. Clearly, $\{1/\nu(K_r)\} \in \ell_2$. With the above settings, Corollary 3.6 yields the results when the suprema over $T$ is taken over the natural numbers $r \geq 1$. Now, for any real $T \geq 1$, with $r = [T]$ the integral part of $T$, we have

$$\frac{\max_{\mathbf{t}\in[-T,T]^s} |\sum_{k=n+1}^m X_k(\omega) e^{i\langle \boldsymbol{\lambda}_k, \mathbf{t}\rangle}|^2}{\log[C_s T^s \cdot m^s \cdot \prod_{i=1}^s (|\lambda_m^{(i)}|^* + 1) + 1] \cdot \sum_{k=n+1}^m |X_k(\omega)|^2}$$

$$\leq 2 \frac{\max_{\mathbf{t}\in[-r-1,r+1]^s} |\sum_{k=n+1}^m X_k(\omega) e^{i\langle \boldsymbol{\lambda}_k, \mathbf{t}\rangle}|^2}{\log[C_s (r+1)^s \cdot m^s \cdot \prod_{i=1}^s (|\lambda_m^{(i)}|^* + 1) + 1] \cdot \sum_{k=n+1}^m |X_k(\omega)|^2}. \quad \square$$

REMARKS. 1. Theorem 7 in [37] is a one-dimensional version of Corollary 3.7 in the periodic case. More specifically, take $s = 1$, and $\{\lambda_n\}$ a nondecreasing sequence of *integers*. Of course, in that case the suprema over $T \geq 1$ is redundant. Recently, Weber proved Theorem 7 without monotonicity, and under the condition $\lambda_n^* \geq n^\delta$, $\delta > 0$ (see also remark 2 after Theorem 3.5)—personal communication. Weber's proof is completely different; it is based on the Dudley metric entropy method, and uses the Borell–Sudakov–Tsirelson inequality (see [22], Lemma 3.1, page 57).

2. By using the reduction principle, one can deduce from Theorem 6 in [16] a weak one-dimensional version of Corollary 3.7. It gives usual integrability (i.e., without the exponential of the square), without the suprema over $T$. The proof of [16] uses a general Gaussian inequality of [18]. Their results hold for centered independent random variables under a quite restrictive condition.

3. As we will see in Section 5, it seems that in applications to random ergodic theory the suprema over $T \geq 1$ is important.

The theorem below shows that the general centered case holds modulo some slight modifications which have no impact on our applications.

NOTATION. For any real numbers $a$ and $b$ we put $a \vee b = \max\{a, b\}$.

THEOREM 3.8. *Let $1 < p \leq 2$, and let $\{X_n\} \subset L_p(\Omega, \mu)$ be a sequence of complex-valued centered independent random variables. Let $\boldsymbol{\lambda}_n = (\lambda_n^{(1)}, \ldots, \lambda_n^{(s)})$* be a sequence of vectors in $\mathbf{R}^s$. Then there exist some positive constants $\varepsilon$ and $C$, independent of $\{X_n\}$, such that (with $0/0$ interpreted as 1)

$$\left\| \sup_{m \geq 1} \max_{0 \leq n < m} \sup_{T \geq 1} \exp\left\{ \varepsilon \cdot \frac{\max_{\mathbf{t} \in [-T,T]^s} |\sum_{k=n+1}^{m} X_k e^{i\langle \boldsymbol{\lambda}_k, \mathbf{t}\rangle}|^2}{\log(T^s \cdot \gamma_m + 1) \cdot (\sum_{k=n+1}^{m} |X_k|^p + \|X_k\|_p^p)^{2/p}} \right\} \right\|_1 \leq C,$$

where $\gamma_m = (12s)^s \cdot m^s \cdot \prod_{i=1}^{s}(|\lambda_m^{(i)}|^* + 1)$. In particular, for a.e. $\omega \in \Omega$ the quantity

$$\sup_{m > n} \sup_{T \geq 1} \left\{ \left( \max_{\mathbf{t} \in [-T,T]^s} \left| \sum_{k=n+1}^{m} X_k(\omega) e^{i\langle \boldsymbol{\lambda}_k, \mathbf{t}\rangle} \right| \right) \right.$$

$$\times \left( \sqrt{\log(T+1) \log\left[m \vee \max_{1 \leq i \leq s}(|\lambda_m^{(i)}|^* + 1) + 1\right]} \right.$$

$$\left. \left. \times \left( \sum_{k=n+1}^{m} |X_k(\omega)|^p + \|X_k\|_p^p \right)^{1/p} \right)^{-1} \right\}$$

is finite.

PROOF. Let $\{X'_n\} \subset L_p(\Omega', \mu')$ be a sequence of independent copies of $\{X_n\}$. Clearly the sequence $\{X_n - X'_n\}$ is symmetric independent on $(\Omega \times \Omega', \mu \otimes \mu')$, and hence Corollary 3.7 applies. So, for $\mu \otimes \mu'$ a.e. $(\omega, \omega') \in \Omega \times \Omega'$, there exists a positive constant $C(\omega, \omega')$, such that for any $m > n \geq 0$, and any $T \geq 1$, we have

(1)
$$\max_{\mathbf{t} \in [-T,T]^s} \left| \sum_{k=n+1}^{m} (X_k(\omega) - X'_k(\omega')) e^{i\langle \boldsymbol{\lambda}_k, \mathbf{t}\rangle} \right|$$
$$\leq C(\omega, \omega') \sqrt{\log(T^s \cdot \gamma_m + 1)} \left( \sum_{k=n+1}^{m} |X_k(\omega) - X'_k(\omega')|^2 \right)^{1/2}.$$

Furthermore, for some universal constants $\varepsilon > 0$ and $C'$, independent of $\{X_n\}$, we have

$$\int \exp\{\varepsilon \cdot [C(\omega, \omega')]^2\} \mu(d\omega) \mu'(d\omega') \leq C'.$$

Let $\mathbf{E}$ and $\mathbf{E}'$ be the corresponding expectations in $\Omega$ and $\Omega'$, respectively. By taking the expectation $\mathbf{E}'$ on the left-hand side of (1), we obtain by Jensen's inequality and $\mathbf{E}'(X'_k) = 0$, that for a.e. $\omega \in \Omega$ we have

$$\mathbf{E}' \left\{ \max_{\mathbf{t} \in [-T,T]^s} \left| \sum_{k=n+1}^{m} (X_k(\omega) - X'_k) e^{i\langle \boldsymbol{\lambda}_k, \mathbf{t}\rangle} \right| \right\}$$

$$\geq \max_{\mathbf{t}\in[-T,T]^s} \mathbf{E}'\left|\sum_{k=n+1}^{m}(X_k(\omega)-X'_k)e^{i\langle\boldsymbol{\lambda}_k,\mathbf{t}\rangle}\right|$$

$$\geq \max_{\mathbf{t}\in[-T,T]^s}\left|\sum_{k=n+1}^{m}X_k(\omega)e^{i\langle\boldsymbol{\lambda}_k,\mathbf{t}\rangle}\right|.$$

By taking the expectation $\mathbf{E}'$ on the right-hand side of (1), and using:
(i) $\|\cdot\|_{\ell_2} \leq \|\cdot\|_{\ell_p}$; (ii) Hölder's inequality; (iii) $(|a|+|b|)^p \leq 2^{p-1}(|a|^p+|b|^p)$, we obtain

$$\mathbf{E}'\left\{C(\omega,\cdot)\sqrt{\log(T^s\cdot\gamma_m+1)}\left(\sum_{k=n+1}^{m}|X_k(\omega)-X'_k|^2\right)^{1/2}\right\}$$

$$\leq \sqrt{\log(T^s\cdot\gamma_m+1)}\,\mathbf{E}'\left\{C(\omega,\cdot)\left(\sum_{k=n+1}^{m}|X_k(\omega)-X'_k|^p\right)^{1/p}\right\}$$

$$\leq \sqrt{\log(T^s\cdot\gamma_m+1)}\,(\mathbf{E}'[C(\omega,\cdot)]^{p/(p-1)})^{(p-1)/p}$$

$$\times \left(\sum_{k=n+1}^{m}\mathbf{E}'|X_k(\omega)-X'_k|^p\right)^{1/p}$$

$$\leq \sqrt{\log(T^s\cdot\gamma_m+1)}\,(\mathbf{E}'[C(\omega,\cdot)]^{p/(p-1)})^{(p-1)/p}\cdot 2^{(p-1)/p}$$

$$\times \left(\sum_{k=n+1}^{m}|X_k(\omega)|^p + \mathbf{E}'|X'_k|^p\right)^{1/p}.$$

By combining that with the previous computation we obtain

(7)
$$\sup_{m>n}\sup_{T\geq 1}\exp\left\{\varepsilon\cdot\left(\max_{\mathbf{t}\in[-T,T]^s}\left|\sum_{k=n+1}^{m}X_k(\omega)e^{i\langle\boldsymbol{\lambda}_k,\mathbf{t}\rangle}\right|^2\right)\right.$$
$$\times\left(2^{2(p-1)/p}\log(T^s\cdot\gamma_m+1)\right.$$
$$\left.\left.\times\left(\sum_{k=n+1}^{m}|X_k(\omega)|^p+\mathbf{E}|X_k|^p\right)^{2/p}\right)^{-1}\right\}$$
$$\leq \exp\{\varepsilon\cdot(\mathbf{E}'[C(\omega,\cdot)]^{p/(p-1)})^{2(p-1)/p}\}.$$

Define the function $\phi_p(u)=\exp\left(u^{2(p-1)/p}\right)$, which is convex in the interval $[K_p,\infty)$, where $K_p:=(\frac{2-p}{2(p-1)})^{p/2(p-1)}$. Using Jensen's inequality and using $(a+b)^\alpha \leq a^\alpha+b^\alpha$, for $a,b\geq 0$ and $0\leq\alpha\leq 1$, we have

$$\exp\{\varepsilon\cdot(\mathbf{E}'[C(\omega,\cdot)]^{p/(p-1)})^{2(p-1)/p}\}$$

$$\begin{aligned}
&= \phi_p(\mathbf{E}'|\varepsilon^{1/2}C(\omega,\cdot)|^{p/(p-1)}) \\
&\leq \phi_p(\mathbf{E}'[|\varepsilon^{1/2}C(\omega,\cdot)|^{p/(p-1)} + K_p]) \\
&\leq \mathbf{E}'[\phi_p(|\varepsilon^{1/2}\cdot C(\omega,\cdot)|^{p/(p-1)} + K_p)] \\
&= \mathbf{E}'\exp\{(|\varepsilon^{1/2}\cdot C(\omega,\cdot)|^{p/(p-1)} + K_p)^{2(p-1)/p}\} \\
&\leq \mathbf{E}'\exp\{\varepsilon\cdot[C(\omega,\cdot)]^2 + K_p^{2(p-1)/p}\} \\
&= \exp\left\{\frac{2-p}{2(p-1)}\right\}\cdot \mathbf{E}'\exp\{\varepsilon\cdot[C(\omega,\cdot)]^2\}.
\end{aligned} \tag{8}$$

By taking the expectation $\mathbf{E}$ in (8) and using the fact that $\mathbf{EE}'\exp\{\varepsilon\times [C(\cdot,\cdot)]^2\} \leq C'$, we obtain the first result with $C = \exp\{\frac{2-p}{2(p-1)}\}\cdot C'$ and by changing the value of $\varepsilon$ to $\varepsilon/2^{2(p-1)/p}$.

The second result follows from the first assertion using the inequality

$$m^s \cdot \prod_{i=1}^{s}(|\lambda_m^{(i)}|^* + 1) \leq \left[m \vee \max_{1\leq i\leq s}(|\lambda_m^{(i)}|^* + 1)\right]^{2s}. \qquad \square$$

**4. Convergence results.** In this section we give general convergence results for sequences of Banach-valued random variables. Then we apply these results with our previous uniform estimates, and we obtain a.e. uniform convergence of random series of almost periodic functions.

Let $\mathbf{B}$ be a separable Banach space with norm $\|\cdot\|$. Let $(\Omega,\mu)$ be a measure space, and let $\mathbf{E}$ be the corresponding expectation. Let $X$ be a random variable on $\Omega$ with values in $\mathbf{B}$. The separability assumption of $\mathbf{B}$ is made in order to avoid measurability complications. For $1\leq p<\infty$, we denote by $\|X\|_p$ the quantity $(\mathbf{E}\|X\|^p)^{1/p}$. The Banach space of all random variables with finite $\|\cdot\|_p$-norm is denoted by $L_p(\Omega,\mu;\mathbf{B})$ or simply by $L_p(\mu;\mathbf{B})$ [or $L_p(\mathbf{B})$]. When $\mathbf{B} = \mathbf{C}$ (or $\mathbf{R}$) we just write $L_p(\Omega,\mu)$ or simply $L_p$.

Let $K$ be a compact metric separable space. In our applications we will take $\mathbf{B} = C(K)$, the space of continuous functions on $K$, with the norm $\|f\| = \max_{x\in K}|f(x)|$ for $f \in C(K)$.

THEOREM 4.1. *Let $\{X_n\} \subset L_p(\Omega,\mu;\mathbf{B})$, with $1<p<\infty$. Let $\{\alpha_n\}$ be a sequence of nonnegative numbers. Let $\gamma$ and $C$ be positive constants, and assume there exists a positive nondecreasing (possibly constant) sequence $\{A_n\}$, with $A_n \leq Cn^\gamma$, such that for every $m > n \geq 0$ we have*

$$\mathbf{E}\left[\left\|\sum_{k=n+1}^{m} X_k\right\|^p\right] \leq A_m \sum_{k=n+1}^{m} \alpha_k. \tag{7}$$

If $\sum_{n=1}^\infty \alpha_n A_n (\log n)^p$ converges, then the series $\sum_{n=1}^\infty X_n$ converges almost everywhere and in $L_p(\Omega, \mu; \mathbf{B})$. Furthermore, we have

$$\left\| \sup_{n \geq 1} \left\| \sum_{k=1}^n X_k \right\| \right\|_p \leq 2e^{1/p}[1 + p^{(p-1)/p}(p+\gamma)] \left( \sum_{n=1}^\infty \alpha_n A_n (\log n)^p \right)^{1/p}.$$

PROOF. By remark 4 after Proposition 2.6, we have

$$\left\| \max_{n < l \leq m} \left\| \sum_{k=n+1}^l X_k \right\| \right\|_p^p \leq A_m (\log_2 4m)^p \sum_{k=n+1}^m \alpha_k,$$

for every $m > n \geq 0$. If we restrict ourselves to $m \geq 2$, by changing $\{A_n\}$, we may replace the $\log_2 4m$ by the natural logarithm $\log m$.

We define a sequence of integers $\{\kappa_n\}$ as the following. Let $\kappa_1$ be the first integer for which $A_{\kappa_1}(\log \kappa_1)^p \geq e$. For every $n \geq 2$, we define inductively

$$\kappa_{n+1} = \max\{m \geq \kappa_n + 1 : A_m(\log m)^p \leq eA_{\kappa_n+1}(\log(\kappa_n + 1))^p\}.$$

Clearly, the sequence $\{\kappa_n\}$ is strictly increasing, and for every $n \geq 1$ we have the following properties:

(i) $A_{\kappa_{n+1}}(\log \kappa_{n+1})^p \leq eA_{\kappa_n+1}(\log(\kappa_n + 1)^p < A_{\kappa_{n+1}+1}(\log(\kappa_{n+1} + 1))^p$;
using (i) we have
(ii) $A_{\kappa_{n+1}} \leq eA_{\kappa_n+1}$ and $A_{\kappa_n+1}(\log(\kappa_n + 1))^p \geq e^n$; by the assumption $A_n \leq Cn^\gamma$ and (ii), we have
(iii) $(p+\gamma)\log(\kappa_n + 1) \geq n$.

(a) Using (7), (ii) and (iii), we obtain

$$\int \sum_{v=1}^\infty v^p \left\| \sum_{k=\kappa_v+1}^{\kappa_{v+1}} X_k \right\|^p d\mu \leq \sum_{v=1}^\infty v^p A_{\kappa_{v+1}} \sum_{k=\kappa_v+1}^{\kappa_{v+1}} \alpha_k$$

$$\leq e(p+\gamma)^p \sum_{n=1}^\infty \alpha_n A_n (\log n)^p < \infty.$$

Hence by Beppo Levi the integrand $\sum_{v=1}^\infty v^p \| \sum_{k=\kappa_v+1}^{\kappa_{v+1}} X_k \|^p$ converges almost everywhere.

(b) For any naturals $r$ and $m$ we obtain, using Hölder's inequality,

$$\left\| \sum_{k=\kappa_m+1}^{\kappa_{m+r}} X_k \right\|^p \leq \left( \sum_{v=m}^{m+r-1} \frac{1}{v} v \left\| \sum_{k=\kappa_v+1}^{\kappa_{v+1}} X_k \right\| \right)^p$$

$$\leq \left( \sum_{v=m}^\infty \frac{1}{v^{p/(p-1)}} \right)^{p-1} \sum_{v=1}^\infty v^p \left\| \sum_{k=\kappa_v+1}^{\kappa_{v+1}} X_k \right\|^p.$$

The first factor in the right-hand side converges to zero as $m \to \infty$, while the last factor converges a.e. by (a), so $\{\sum_{k=1}^{\kappa_m} X_k\}$ is a Cauchy sequence a.e., and hence converges a.e. By taking integrals of the above inequality, and considering the convergence proved in (a), $\{\sum_{k=1}^{\kappa_m} X_k\}$ is a Cauchy sequence in $L_p(\Omega, \mu; \mathbf{B})$-norm, and hence converges in norm.

(c) Using (7) and (i), we have

$$\sum_{m=1}^{\infty} \int \max_{\kappa_m < n \leq \kappa_{m+1}} \left\| \sum_{k=\kappa_m+1}^{n} X_k \right\|^p d\mu \leq \sum_{m=1}^{\infty} A_{\kappa_{m+1}} (\log \kappa_{m+1})^p \sum_{k=\kappa_m+1}^{\kappa_{m+1}} \alpha_k$$

$$\leq e \sum_{n=1}^{\infty} \alpha_n A_n (\log n)^p < \infty.$$

Now, (b) and (c) imply that $\sum_{n=1}^{\infty} X_n$ converges a.e. in $\mathbf{B}$ to $X := \lim_{m \to \infty} \sum_{n=1}^{\kappa_m} X_n$, since for $\kappa_m < n \leq \kappa_{m+1}$, we have

$$\left\| \sum_{k=1}^{n} X_k - X \right\| \leq \left\| \sum_{k=1}^{\kappa_m} X_k - X \right\| + \left\| \sum_{k=\kappa_m+1}^{n} X_k \right\|$$

$$\leq \left\| \sum_{k=1}^{\kappa_m} X_k - X \right\| + \max_{\kappa_m < n \leq \kappa_{m+1}} \left\| \sum_{k=\kappa_m+1}^{n} X_k \right\|.$$

By considering the norm convergence proved in (a) and (b), the $L_p(\Omega, \mu; \mathbf{B})$-norm convergence follows by taking the $L_p(\Omega, \mu)$-norm in the above inequality.

Now we will prove that $\sup_{n \geq 1} \| \sum_{k=1}^{n} X_k \| \in L_p(\mu)$. The inequality in (b) with $m = 1$ yields

$$\sup_{r \geq 1} \left\| \sum_{k=\kappa_1+1}^{\kappa_{r+1}} X_k \right\|^p \leq \left( \sum_{v=1}^{\infty} \frac{1}{v^{p/(p-1)}} \right)^{p-1} \sum_{v=1}^{\infty} v^p \left\| \sum_{k=\kappa_v+1}^{\kappa_{v+1}} X_k \right\|^p.$$

Integration of the above inequality and application of (a) yield

(1) $$\int \sup_{r \geq 1} \left\| \sum_{k=\kappa_1+1}^{\kappa_{r+1}} X_k \right\|^p d\mu \leq p^{p-1} e(p+\gamma)^p \sum_{n=1}^{\infty} \alpha_n A_n (\log n)^p < \infty.$$

The inequality in (c) yields

(7) $$\int \sup_{m \geq 1} \max_{\kappa_m < n \leq \kappa_{m+1}} \left\| \sum_{k=\kappa_m+1}^{n} X_k \right\|^p d\mu \leq \int \sum_{m=1}^{\infty} \max_{\kappa_m < n \leq \kappa_{m+1}} \left\| \sum_{k=\kappa_m+1}^{n} X_k \right\|^p d\mu$$

$$\leq e \sum_{n=1}^{\infty} \alpha_n A_n (\log n)^p < \infty.$$

Using (1) and (7), we have

$$\left(\int \sup_{n>\kappa_1}\left\|\sum_{k=\kappa_1+1}^n X_k\right\|^p d\mu\right)^{1/p}$$
$$\leq e^{1/p}[1+p^{(p-1)/p}(p+\gamma)]\left(\sum_{n=1}^\infty \alpha_n A_n(\log n)^p\right)^{1/p}.$$

Since a.e. we have

$$\sup_{n\geq 1}\left\|\sum_{k=1}^n X_k\right\| \leq \|X_1\| + \|X_2\| + \cdots + \|X_{\kappa_1}\| + \sup_{n>\kappa_1}\left\|\sum_{k=\kappa_1+1}^n X_k\right\|,$$

we have

$$\left(\int \sup_{n\geq 1}\left\|\sum_{k=1}^n X_k\right\|^p d\mu\right)^{1/p}$$
$$\leq A_{\kappa_1}(\alpha_1 + \cdots + \alpha_{\kappa_1})$$
$$+ e^{1/p}[1+p^{(p-1)/p}(p+\gamma)]\left(\sum_{n=1}^\infty \alpha_n A_n(\log n)^p\right)^{1/p}.$$

So, we obtain the maximal inequality. □

REMARKS. 1. Theorem 4.1 uses techniques of Theorem 2.4 in [8]. In [8] the theorem was obtained without the condition $A_n \leq Cn^\gamma$ (which is not restrictive in our applications), but under the assumption $\sum_{n=1}^\infty A_{2^n}\alpha_n(\log n)^p$.

2. If (7) holds for $\{A_m\}$ bounded, then the condition $\sum_{n=1}^\infty \alpha_n < \infty$ is not sufficient in general for the a.e. convergence (see Menchoff's example [25], Theorem 3).

THEOREM 4.2. *Let $\{X_n\} \subset L_p(\Omega, \mu; \mathbf{B})$, with $1 < p < \infty$. Let $\{\alpha_n\}$ be a sequence of nonnegative numbers. Assume there exists a positive nondecreasing unbounded sequence $\{A_n\}$, with $A_1 \geq 1$, and some positive constant $q \geq 1$, such that for every $m > n \geq 0$ we have*

$$(8) \quad \mathbf{E}\left[\max_{n<l\leq m}\left\|\sum_{k=n+1}^l X_k\right\|^p\right] \leq A_m\left(\sum_{k=n+1}^m \alpha_k\right)^q.$$

*If the series*

$$\sum_{n=1}^\infty 2^{n/p}\left(\sum_{\{k:\, 2^n\leq A_k\leq 2^{n+1}\}} \alpha_k\right)^{q/p},$$

*or in particular*

$$\sum_{n=2}^{\infty} \frac{(\sum_{\{k:\, 2^{A_k} \geq n\}} \alpha_k)^{q/p}}{n(\log n)^{1-1/p}}$$

*converges, then the series* $\sum_{n=1}^{\infty} X_n$ *converges almost everywhere and in* $L_p(\Omega, \mu; \mathbf{B})$. *Furthermore, we have*

$$\left\| \sup_{n \geq 1} \left\| \sum_{k=1}^{n} X_k \right\| \right\|_p \leq 2 \sum_{n=1}^{\infty} 2^{n/p} \left( \sum_{\{k:\, A_k \geq 2^n\}} \alpha_k \right)^{s/p}.$$

PROOF. Define two sequences $\{\kappa_n\}$ and $\{l_n\}$ by induction. Let $\kappa_1 = l_1 = 0$. Let $l_{n+1}$ be the integer such that $2^{l_{n+1}} \leq A_{\kappa_n+1} < 2^{l_{n+1}+1}$. Then define $\kappa_{n+1} = \max\{m \geq \kappa_n + 1 : 2^{l_{n+1}} \leq A_m < 2^{l_{n+1}+1}\}$. Clearly, $\kappa_n$ is finite by the unboundness and nondecreasingness of $\{A_n\}$; moreover, it is increasing. In particular, by assumption on $A_1$, we have $l_2 = 0$. Now, for every $n \geq 1$ we have

$$\left( \int \max_{\kappa_n < m \leq \kappa_{n+1}} \left\| \sum_{k=\kappa_n+1}^{m} X_k \right\|^p d\mu \right)^{1/p} \leq A_{\kappa_{n+1}}^{1/p} \left( \sum_{k=\kappa_n+1}^{\kappa_{n+1}} \alpha_k \right)^{q/p}$$

$$\leq 2^{(l_{n+1}+1)/p} \left( \sum_{2^{l_{n+1}} \leq A_k \leq 2^{l_{n+1}+1}} \alpha_k \right)^{q/p}.$$

Hence,

$$\sum_{n=1}^{\infty} \left\| \max_{\kappa_n < m \leq \kappa_{n+1}} \left\| \sum_{k=\kappa_n+1}^{m} X_k \right\| \right\|_p \leq 2 \sum_{n=0}^{\infty} 2^{n/p} \left( \sum_{2^n \leq A_k \leq 2^{n+1}} \alpha_k \right)^{q/p}.$$

This proves the a.e. convergence.

Let $\kappa_n < m \leq \kappa_{n+1}$. We have

$$\left\| \sum_{k=1}^{m} X_k \right\| \leq \sum_{l=1}^{n-1} \left\| \sum_{\kappa_l+1}^{\kappa_{l+1}} X_k \right\| + \max_{\kappa_n+1 \leq l \leq \kappa_{n+1}} \left\| \sum_{k=\kappa_n+1}^{l} X_k \right\|.$$

The maximal inequality then follows from the previous computations. □

REMARKS. 1. For $q > 1$, using remarks 1 and 2 after Proposition 2.6, we could assume the following weaker condition instead of (8):

(9) $$\left\| \sum_{k=n+1}^{m} X_k \right\|_p^p \leq A_m \left( \sum_{k=n+1}^{m} \alpha_k \right)^q \quad \text{for every } m > n \geq 0.$$

2. Clearly, $\mu(\sup_{j \geq 1} \|\sum_{k=n+1}^{n+j} X_k\| \geq \varepsilon) \xrightarrow{n} 0$ is equivalent to the $\mu$-a.e. convergence of $\sum_{n=1}^{\infty} X_n$. Hence, if for some $1 \leq p < \infty$ and $1 \leq q < \infty$, (8) holds for every $m > n \geq 0$ with $\{A_n\}$ bounded, then the condition $\sum_{n=1}^{\infty} \alpha_n < \infty$ is sufficient for the convergence result of the above theorem. Moreover, we have

$$\mathbf{E}\left[\sup_{n\geq 1}\left\|\sum_{k=1}^{n} X_k\right\|^p\right] \leq \sup_{n\geq 1}\{A_n\}\left(\sum_{k=1}^{\infty} \alpha_k\right)^q. \tag{10}$$

3. Using remarks 1 and 2 above, the condition $\sum_{n=1}^{\infty} \alpha_n < \infty$ is sufficient for the convergence result of the above theorem under (9), with $q > 1$ and $\{A_m\}$ bounded. In that case (10) holds by multiplying the right-hand side by factor $(1 - \frac{1}{2^{(q-1)/p}})^{-p}$ (see also [8], Theorem 2.5).

4.1. *Convergence of random almost periodic series.* Let $\{X_n\} \subset L_\infty(\Omega, \mu)$ be a sequence of complex centered random variables. In the following we will use the notation:

$$\alpha_i = \sum_{k=1}^{i} \|X_k \mathbf{E}(X_i|\mathcal{F}_k)\|_\infty \qquad \text{for every } i \geq 1.$$

THEOREM 4.3. *Let $\{X_n\} \subset L_\infty(\Omega, \mu)$ be a sequence of complex centered random variables, and let $2 \leq p < \infty$. Let $\{f_n\}$ be a $\{\sigma_n\}$-system on $(K, \nu)$. If the series*

$$\sum_{n=2}^{\infty} \frac{(\sum_{\{k \,:\, 2^{(\log(\sigma_k+1))^{p/2}} \geq n\}} \alpha_k)^{1/2}}{n(\log n)^{1-1/p}}$$

*converges, then for a.e. $\omega \in \Omega$ the series $\sum_{n=1}^{\infty} X_n(\omega) f_k$ converges uniformly on $K$. Furthermore, it converges in $L_p(\Omega, \mu; C(K))$, and $\sup_{n\geq 1} \max_{x\in K} |\sum_{k=1}^{n} X_k f_k(x)|$ is in $L_p(\mu)$.*

PROOF. Put $\mathbf{B} = C(K)$. Using the second part of Theorem 3.1, up to appropriate constant, we obtain (8) with $p \geq 2$, $q = p/2 \geq 1$, $A_n = C_p(\log(\sigma_n + 1))^{p/2}$ and $\{\alpha_i\}$ as defined above. Theorem 4.2 yields the results. □

THEOREM 4.4. *Let $\{X_n\} \subset L_\infty(\Omega, \mu)$ be a martingale difference sequence. Let $\boldsymbol{\lambda}_n = (\lambda_n^{(1)}, \ldots, \lambda_n^{(s)})$ be a sequence in $\mathbf{R}^s$, and put $\gamma_k = k \vee \max_{1 \leq i \leq s} |\lambda_k^{(i)}|^*$. If the series*

$$\sum_{n=2}^{\infty} \frac{(\sum_{\{k \,:\, \gamma_k \geq n\}} \|X_k\|_\infty^2)^{1/2}}{n(\log n)^{1/2}}$$

*converges, then for every $T \geq 1$ and for almost every $\omega \in \Omega$, the series $\sum_{n=1}^{\infty} X_n(\omega) e^{i\langle \boldsymbol{\lambda}_n, \mathbf{t}\rangle}$ converges uniformly in $\mathbf{t} \in [-T,T]^s$. Furthermore, it converges in $L_2(\Omega, \mu; C([-T,T]^s))$, and*

$$\sup_{n \geq 1} \max_{\mathbf{t} \in [-T,T]^s} \left| \sum_{k=1}^{n} X_k e^{i\langle \boldsymbol{\lambda}_n, \mathbf{t}\rangle} \right|$$

*is in $L_2(\mu)$.*

PROOF. Let $K = [-T,T]^s$, and let $\nu$ be the Lebesgue measure on $\mathbf{R}^s$. By Example 2.1 the sequence $\{e^{i\langle \boldsymbol{\lambda}_n, \mathbf{t}\rangle}\}$ forms a $\{(n \vee \max_{1 \leq i \leq s} |\lambda_n^{(i)}|^*)^{2s}\}$-system on $(K, \nu)$. With the above settings, Theorem 4.3 with $p=2$ yields the assertions of the theorem if

$$\sum_{n=2}^{\infty} \frac{(\sum_{\{k\,:\,\gamma_k \geq n^{1/2s}\}} \|X_k\|_{\infty}^2)^{1/2}}{n(\log n)^{1/2}} < \infty.$$

By a change of variable in the series above we obtain the condition of the theorem. □

THEOREM 4.5. *Let $\{X_n\} \subset L_2(\Omega, \mu)$ be a sequence of centered independent random variables. Let $\boldsymbol{\lambda}_n = (\lambda_n^{(1)}, \ldots, \lambda_n^{(s)})$ be a sequence of vectors in $\mathbf{R}^s$, and put $\gamma_k = k \vee \max_{1 \leq i \leq s} |\lambda_k^{(i)}|^*$. If the series*

$$(11) \qquad \sum_{n=2}^{\infty} \frac{(\sum_{\{k\,:\,\gamma_k \geq n\}} \|X_k\|_2^2)^{1/2}}{n(\log n)^{1/2}}$$

*converges, then for every $T \geq 1$ and for almost every $\omega \in \Omega$, the series $\sum_{n=1}^{\infty} X_n(\omega) e^{i\langle \boldsymbol{\lambda}_n, \mathbf{t}\rangle}$ converges uniformly in $\mathbf{t} \in [-T,T]^s$. Furthermore, it converges in $L_2(\Omega, \mu; C([-T,T]^s))$ and*

$$\sup_{n \geq 1} \max_{\mathbf{t} \in [-T,T]^s} \left| \sum_{k=1}^{n} X_k e^{i\langle \boldsymbol{\lambda}_n, \mathbf{t}\rangle} \right|$$

*is in $L_2(\mu)$.*

PROOF. By Example 2.1 the sequence $\{e^{i\langle \boldsymbol{\lambda}_n, \mathbf{t}\rangle}\}$ forms a $\{(n \vee \max_{1 \leq i \leq s} |\lambda_n^{(i)}|^*)^{2s}\}$-system on $([-T,T], \nu)$. Using Colorally 3.3 we obtain inequality (8) with $p=2$, $q=1$ and $A_n = 2sC_p \log \gamma_n$. Theorem 4.2 yields the result, as in Theorem 4.4. □

REMARKS. 1. The a.e. uniform convergence of multidimensional almost periodic series with square integrable *symmetric* independent coefficients was considered in Marcus and Pisier [23], Chapter VII, Section 1. Their proofs use the metric entropy method. Using the metric entropy method one can obtain a more precise condition than (11). Our condition is the same sufficient condition as it is implicit in [23], Chapter VII, Lemma 1.1. In any case, Theorem 4.5 completely recovers Theorem 5.1.5 of [34].

2. In [24] the a.e. uniform convergence of multidimensional almost periodic series with i.i.d. $p$-stable coefficients was considered. As it is implicit in [24], Remark 4.4 (use [23], Chapter VII, Lemma 1.1 and then [24], Theorem B) the result is valid for symmetric i.i.d. coefficients [not necessarily $p$-stable or in $L_p(\mu)$].

3. Extending the results beyond the scope of symmetric random variable is possible, under various conditions, by a symmetrization procedure (see, e.g., [5]).

4. Let $1 < p < 2$, and let $\{X_n\} \subset L_p(\mu)$ be centered independent random variables. We could use our previous results in order to obtain a.e. uniform convergence also in this case (i.e., for $p < 2$). In that case our derived sufficient condition is not as good as given in [24], Remark 4.4.

**5. Applications to ergodic theory.** Let $(Y, \Sigma, \pi)$ be a probability space, and let $V_1, \ldots, V_s$ be pairwise commuting isometries on $L_2(Y, \pi)$. For every $\mathbf{j} := (j^{(1)}, \ldots, j^{(s)}) \in \mathbf{N}^s$ and $f \in \mathcal{H} := L_2(Y, \pi)$, define the action $V^{\mathbf{j}} f = V_1^{j^{(1)}} \circ \cdots \circ V_s^{j^{(s)}} f$. Using the dilation theorem for pairwise commuting family of isometries (see, e.g., [28], Chapter I, Proposition 6.2), there exist a Hilbert space $\mathcal{H}' \supset \mathcal{H}$ and a family $U_1, \ldots, U_s$ of pairwise commuting *unitary* operators on $\mathcal{H}'$, such that for every $\mathbf{j} \in \mathbf{N}^s$ and $f \in \mathcal{H}$ we have $V^{\mathbf{j}} f = P_{\mathcal{H}} U^{\mathbf{j}} f$, where $P_{\mathcal{H}}$ is the orthogonal projection of $\mathcal{H}'$ onto $\mathcal{H}$ and $U^{\mathbf{j}} f = U_1^{j^{(1)}} \circ \cdots \circ U_s^{j^{(s)}} f$. By the spectral theorem for the unitary representation $\{U^{\mathbf{k}} : \mathbf{k} \in \mathbf{Z}^s\}$ (e.g., see [31], Chapter X, Section 140), there is a positive measure $\mu_f$ on $[-\pi, \pi)^s$, called the spectral measure of $f$, such that for any $\mathbf{k} \in \mathbf{Z}^s$ we have

$$\langle U^{\mathbf{k}} f, f \rangle = \int_{[-\pi,\pi)^s} e^{i \langle \mathbf{k}, \mathbf{t} \rangle} \, d\mu_f(\mathbf{t}),$$

where $\langle \mathbf{k}, \mathbf{t} \rangle$ denotes the inner product in $\mathbf{R}^s$.

Let $\{\alpha_k\}$ be a sequence of complex numbers, and let $\{\mathbf{j}_k\} \subset \mathbf{N}^s$. Using the dilation theorem and the spectral theorem, for every $n$ we have

$$(12) \quad \left\| \sum_{k=1}^n \alpha_k V^{\mathbf{j}_k} f \right\|_2 = \left\| P_{\mathcal{H}} \left( \sum_{k=1}^n \alpha_k U^{\mathbf{j}_k} f \right) \right\|_2 \leq \|f\|_2 \max_{\mathbf{t} \in [-\pi,\pi)^s} \left| \sum_{k=1}^n \alpha_k e^{i \langle \mathbf{j}_k, \mathbf{t} \rangle} \right|.$$

REMARKS. 1. Let $T_1$ and $T_2$ be commuting contractions on $\mathcal{H}$. One can consider Ando's unitary dilation for pairs of contractions ([28], Chapter I, Theorem 6.4). Specifically, there exist two commuting unitary operators $U_1$ and $U_2$, acting on $\mathcal{H}' \supset \mathcal{H}$, such that for every natural number $n$ and $m$ we have $T_1^n T_2^m f = P_{\mathcal{H}} U_1^n U_2^m f$. In the case of two commuting contractions, (12) is still true when $V_i$ are replaced by $T_i$.

2. Parrott [30] gave an example (see also [28], Chapter I, Section 3) which showed that the dilation theorem is no longer true in the case of more than two commuting contractions (for an analogue for commuting Markov operators see [15]).

3. For more than two commuting contractions, additional conditions should be added in order to obtain a *regular dilation* (see [28], Chapter I, Theorem 9.1 or a theorem of Brehmer [27], Chapter 6). A simple condition that one can assume on a family of contractions is the pairwise doubly permutability (where doubly stands for commuting also with the conjugate); see [28], Chapter I, Proposition 9.2.

NOTATION. For $\mathbf{j} = (j^{(1)}, \ldots, j^{(s)}) \in \mathbf{Z}^d$, put $|\mathbf{j}| = \max\{|j^{(1)}|, \ldots, |j^{(s)}|\}$. We recall our notation $a \vee b = \max\{a, b\}$, and $c_m^* = \max_{1 \leq n \leq m} c_n$ for a positive sequence $\{c_n\}$. If $\mathbf{j}_n = (j_n^{(1)}, \ldots, j_n^{(s)})$ is a sequence of vectors, then by our notation we have

$$|\mathbf{j}_m|^* = \max_{1 \leq i \leq s} |j_m^{(i)}|^* = \max_{1 \leq n \leq m} \max_{1 \leq i \leq s} |j_n^{(i)}|.$$

THEOREM 5.1. *Let $\{X_n\} \subset L_2(\Omega, \mu)$ be a sequence of centered independent random variables, and let $\{\mathbf{j}_n\} \subset \mathbf{N}^s$. If the series $\sum_{n=1}^{\infty} \|X_n\|_2^2 \log(n \vee |\mathbf{j}_n|^*)(\log n)^2$ converges, then there exists a set of full measure $\Omega^* \subset \Omega$, such that for every $\omega \in \Omega^*$, for every commuting family of isometries $V_1, \ldots, V_s$ on a space $L_2(Y, \pi)$ and any $f \in L_2(\pi)$, the series*

$$\sum_{n=1}^{\infty} X_n(\omega) V^{\mathbf{j}_n} f \qquad \text{converges } \pi\text{-a.e.}$$

*Furthermore, for every $\omega \in \Omega^*$ there exists a constant $K_\omega < \infty$, independent of $f$, such that*

$$\left\| \sup_{n \geq 1} \left| \sum_{k=1}^{n} X_k(\omega) V^{\mathbf{j}_k} f \right| \right\|_2 \leq K_\omega \|f\|_2.$$

PROOF. We construct from $\{X_n\}$ two sequences of random variables in the following way:

$$Y_n = \begin{cases} X_n, & \text{if } |\mathbf{j}_n|^* \leq e^n, \\ 0, & \text{otherwise,} \end{cases}$$

and

$$Z_n = \begin{cases} X_n, & \text{if } |\mathbf{j}_n|^* > e^n, \\ 0, & \text{otherwise.} \end{cases}$$

Hence $X_n = Y_n + Z_n$.

Let $\{n_k\}$ be the sequence of integers for which $|\mathbf{j}_{n_k}|^* > e^{n_k}$, that is, the sequence for which $Z_{n_k}$ is not null. Notice that for every $\omega \in \Omega$, we have

$$
\begin{aligned}
\sum_{n=2}^{\infty} |Z_n(\omega)| &\leq \left( \sum_{k=1}^{\infty} |Z_{n_k}(\omega)|^2 \log(n_k \vee |\mathbf{j}_{n_k}|^*)(\log n_k)^2 \right)^{1/2} \\
&\quad \times \left( \sum_{k=1}^{\infty} \frac{1}{\log(n_k \vee |\mathbf{j}_{n_k}|^*)(\log n_k)^2} \right)^{1/2} \\
&\leq \left( \sum_{n=2}^{\infty} |X_n(\omega)|^2 \log(n \vee |\mathbf{j}_n|^*)(\log n)^2 \right)^{1/2} \\
&\quad \times \left( \sum_{n=2}^{\infty} \frac{1}{n(\log n)^2} \right)^{1/2}.
\end{aligned}
\tag{13}
$$

By our assumption and the theorem of Beppo Levi, for $\mu$ a.e. $\omega \in \Omega$ we have

$$
\sum_{n=1}^{\infty} (|X_n(\omega)|^2 + \|X_n\|_2^2) \log(n \vee |\mathbf{j}_n|^*)(\log n)^2 < \infty. \tag{1}
$$

Hence, by (13) and (1), it remains only to consider the series $\sum_{n=1}^{\infty} Y_n V^{\mathbf{j}_n} f$. Since $Y_n$ is null when $|\mathbf{j}_n|^* > e^n$, we may and do assume from now on (modifying $\{\mathbf{j}_n\}$ when necessary) that *for every $n \geq 1$, we have $|\mathbf{j}_n|^* \leq e^n$*.

On the other hand, by the second assertion of Theorem 3.8, for $\mu$-a.e. $\omega \in \Omega$ there exists a constant $C_\omega$, such that for every $m > n \geq 1$,

$$
\sup_{\mathbf{t} \in [-\pi, \pi)^s} \left| \sum_{k=n+1}^{m} Y_k(\omega) e^{i\langle \mathbf{j}_k, \mathbf{t} \rangle} \right|^2
$$
$$
\leq C_\omega \log(m \vee |\mathbf{j}_m|^* + 1) \sum_{k=n+1}^{m} |Y_k(\omega)|^2 + \|Y_k\|_2^2.
\tag{7}
$$

Let $\Omega^*$ be the set of $\omega$ for which (1) and (7) hold, and fix $\omega \in \Omega^*$.

Using (7) together with the spectral theorem [see (12)], we have for every $m > n \geq 0$,

$$
\left\| \sum_{k=n+1}^{m} Y_k(\omega) V^{\mathbf{j}_k} f \right\|_2^2 \leq \|f\|_2^2 C_\omega \log(m \vee |\mathbf{j}_m|^* + 1) \sum_{k=n+1}^{m} |Y_k(\omega)|^2 + \|Y_k\|_2^2.
$$

Hence, the condition in (7) is satisfied with $\alpha_n = |Y_n(\omega)|^2 + \|Y_n\|_2^2$ and $A_n^{(\omega)} = \|f\|_2^2 C_\omega \log(n \vee |\mathbf{j}_n|^* + 1)$. Since $|\mathbf{j}_n|^* \leq e^n$ we are allowed to use Theorem 4.1. Using (1), Theorem 4.1 yields the two assertions of the theorem for the sequence $\{Y_n\}$. □

THEOREM 5.2. *Let $\{X_n\} \subset L_2(\Omega, \mu)$ be a sequence of centered independent random variables, and let $\{p_n\}$ and $\{q_n\}$ be sequences of natural numbers. If the series*

$$\sum_{n=1}^{\infty} \|X_n\|_2^2 \log(n \vee p_n^* \vee q_n^*)(\log n)^2$$

*converges, then there exists a set of full measure $\Omega^* \subset \Omega$, such that for every $\omega \in \Omega^*$, for every commuting contraction $T_1$ and $T_2$ on a space $L_2(Y, \pi)$ and any $f \in L_2(\pi)$, the series*

$$\sum_{n=1}^{\infty} X_n(\omega) T_1^{p_n} T_2^{q_n} f \qquad \text{converges } \pi\text{-a.e.}$$

*Furthermore, for every $\omega \in \Omega^*$ there exists a constant $K_\omega < \infty$, independent of $f$, such that*

$$\left\| \sup_{n \geq 1} \left| \sum_{k=1}^{n} X_k(\omega) T_1^{p_n} T_2^{q_n} f \right| \right\|_2 \leq K_\omega \|f\|_2.$$

PROOF. As we mentioned in the remarks at the beginning of the section, we can use the unitary dilation in the case of two commuting contractions. In this case also (12) is still true. We proceed as in Theorem 5.1. □

REMARKS. 1. Theorem 5.2 extends the part of [8], Theorem 4.2, related to square integrable $\{X_n\}$, to the case of two commuting contractions. Moreover, the convergence is along certain subsequences.

2. Using Theorems 3.8 and 4.1, we could consider the case $\{X_n\} \subset L_p(\Omega, \mu)$, $1 < p < 2$. It turns out that Theorem 4.1 does not lead to the best result that one can obtain. In addition to Theorem 3.8 another tool seems to be needed. It will be done in the forthcoming paper [6].

Let $\tau_1, \ldots, \tau_s$ be pairwise commuting measure-preserving transformations of a probability space $(Y, \Sigma, \pi)$. For any $\mathbf{j} \in \mathbf{N}^s$ and every $\Sigma$-measurable $f$, we define $T^{\mathbf{j}} f = f \circ \tau_1^{j^{(1)}} \circ \cdots \circ \tau_s^{j^{(s)}}$. By assumptions, for any $\mathbf{j} \in \mathbf{N}^s$ the operator $T^{\mathbf{j}}$ is an isometry of $L_q(Y, \pi)$ for any $1 \leq q < \infty$. Moreover, for fixed $1 \leq q < \infty$ the action $\{T^{\mathbf{j}} : j \in \mathbf{N}^s\}$ on $L_q(Y, \pi)$ is an isometric representation of the semigroup $\mathbf{N}^s$.

As in [33] (and [7]) we want to use Stein's complex interpolation in order to obtain a.e. convergence results, under the assumptions of Theorem 5.1, also for $f \in L_q(Y, \pi)$, $1 < q < 2$.

THEOREM 5.3. *Let $\{X_n\} \subset L_2(\Omega, \mu)$ be a sequence of centered independent random variables, and let $\{\mathbf{j}_n\} \subset \mathbf{N}^s$. If the series $\sum_{n=1}^{\infty} \|X_n\|_2^2 \log(n \vee$* $|\mathbf{j}_n|^*)(\log n)^2$ *converges, then there exists a set of full measure* $\Omega^* \subset \Omega$, *such that for every* $\omega \in \Omega^*$, *for every commuting family of measure-preserving transformations* $\tau_1, \ldots, \tau_s$ *on* $(Y, \pi)$ *and any* $f \in L_q(Y, \pi)$, $1 < q \leq 2$, *the series*

$$\sum_{n=1}^{\infty} \frac{X_n(\omega) T^{\mathbf{j}_n} f}{n^{(2-q)/2q}} \qquad \text{converges } \pi\text{-a.e.}$$

*Furthermore, for every* $\omega \in \Omega^*$ *there exists a constant* $K_\omega < \infty$, *independent of* $f$, *such that*

$$\left\| \sup_{n \geq 1} \left| \sum_{k=1}^{n} \frac{X_k(\omega) T^{\mathbf{j}_k} f}{n^{(2-q)/2q}} \right| \right\|_q \leq K_\omega \|f\|_q.$$

PROOF. As noted in the proof of Theorem 5.1, we may and do assume that $|\mathbf{j}_n|^* \leq e^n$ for every $n \geq 1$.

By the Beppo Levi theorem and by our assumptions, for $\mu$-a.e. $\omega \in \Omega$ we have

(1)	$$\sum_{n=1}^{\infty} (|X_n(\omega)|^2 + \|X_n\|_2^2) \log(n \vee |\mathbf{j}_n|^*)(\log n)^2 < \infty.$$

By the second assertion in Theorem 3.8, for $\mu$-a.e. $\omega \in \Omega$, there exists $C_\omega$ such that, for every $m > n \geq 1$ and every $K \geq 1$, we have

$$\max_{(\eta, \mathbf{t}) \in [-K, K] \times [-K, K]^s} \left| \sum_{k=n+1}^{m} X_k(\omega) e^{-i(1/2)\eta \log k} e^{i \langle \mathbf{j}_k, \mathbf{t} \rangle} \right|^2$$

$$\leq C_\omega \log(K+1) \log(m \vee \tfrac{1}{2}|\mathbf{j}_m|^* \log m + 1) \sum_{k=n+1}^{m} |X_k(\omega)|^2 + \|X_k\|_2^2.$$

In particular, for every $\eta \in \mathbf{R}$ we have

(7)	$$\max_{\mathbf{t} \in [-\pi, \pi)^s} \left| \sum_{k=n+1}^{m} X_k(\omega) e^{-i(1/2)\eta \log k} e^{i \langle \mathbf{j}_k, \mathbf{t} \rangle} \right|^2$$

$$\leq 2 C_\omega \log(|\eta| + \pi) \log(m \vee |\mathbf{j}_m|^* + 1) \sum_{k=n+1}^{m} |X_k(\omega)|^2 + \|X_k\|_2^2.$$

Let $\Omega^*$ be the set for which (1) and (7) hold, and fix $\omega \in \Omega^*$. Hence, for any $f \in L_2(Y, \pi)$ we deduce from (7) and (12), that

(14)	$$\left\| \sum_{k=n+1}^{m} X_k(\omega) e^{-i(1/2)\eta \log k} T^{\mathbf{j}_k} f \right\|_2^2$$

$$\leq 2 C_\omega \|f\|_2^2 \log(|\eta| + \pi) \log(m \vee |\mathbf{j}_m|^* + 1) \sum_{k=n+1}^{m} |X_k(\omega)|^2 + \|X_k\|_2^2.$$

For any complex $\zeta = \xi + i\eta$ with $0 \leq \xi \leq 1$ put $\Psi_{n,\zeta}(T) := \sum_{k=1}^{n} X_k(\omega) \times k^{-(1/2)\zeta} T^{\mathbf{j}_k}$. Using (1) and (14), we apply Theorem 4.1 with $\alpha_n = |X_n(\omega)|^2 + \|X_n\|_2^2$ and $A_n = 2C_\omega \|f\|_2^2 \log(|\eta| + \pi) \log(n \vee |\mathbf{j}_n|^* + 1)$ in order to obtain

$$\left\| \sup_{n \geq 1} |\Psi_{n,i\eta}(T)f| \right\|_2 \leq C_1 \sqrt{\log(|\eta| + \pi)} \|f\|_2,$$

for some $C_1 > 0$, which does not depend on $\eta$ or $f$. On the other hand we have

$$\| \sup_{n \geq 1} |\Psi_{n,1+i\eta}(T)f| \|_1$$

(15)
$$\leq \|f\|_1 \left( \sum_{n=1}^{\infty} |X_n|_2^2 \log(n \vee |\mathbf{j}_n|^*)(\log n)^2 \right)^{1/2}$$

$$\times \left( \sum_{n=1}^{\infty} \frac{1}{n \log(n \vee |\mathbf{j}_n|^*)(\log n)^2} \right)^{1/2}.$$

For a subset $A \subset Y$ let $M_A$ be the operator of multiplication by $\mathbf{1}_A$. For any bounded integer-valued function $I \geq 1$ defined on $Y$ we have linear operators

$$\Psi_{I,\zeta}(T) = \sum_{j=1}^{\max I} M_{\{I=j\}} \sum_{k=1}^{j} X_k(\omega) k^{-(1/2)\zeta} T^{\mathbf{j}_k}.$$

Hence for $f \in L_1(Y, \pi)$ and $y \in Y$ we have

$$\Psi_{I,\zeta}(T)f(y) = \sum_{k=1}^{I(y)} X_k(\omega) k^{-(1/2)\zeta} T^{\mathbf{j}_k} f(y),$$

so $|\Psi_{I,\zeta}(T)f(y)| \leq \sup_{n \geq 1} |\sum_{k=1}^{n} X_k(\omega) k^{-(1/2)\zeta} T^{\mathbf{j}_k} f(y)|$. For $f_1$ and $f_2$ simple functions on $Y$ it is easy to check that $\Phi(\zeta) := \int \Psi_{I,\zeta}(T) f_1 \cdot f_2 \, d\pi$ is continuous in the strip $0 \leq \xi \leq 1$ and analytic in its interior.

Our previous estimates yield

$$\||\Psi_{I,i\eta}(T)f|\|_2 \leq C_1 \|f\|_2 \sqrt{\log(|\eta| + \pi)}, \qquad f \in L_2(\pi),$$

$$\||\Psi_{I,1+i\eta}(T)f|\|_1 \leq C_2 \|f\|_1, \qquad f \in L_1(\pi).$$

Stein's interpolation ([38], Theorem XII.1.39) yields that for $1 \leq q \leq 2$ we have, with $t = (2-q)/q$,

$$\||\Psi_{I,t}(T)f|\|_q \leq C_{1,2,q} \|f\|_q, \qquad f \in L_q(\pi).$$

Note that the constant $C_{1,2,q}$ depends only on $C_1$, $C_2$ and $q$, and not on the choice of $I$ or $f$. Keeping $q$ fixed and taking $f \in L_q(\pi)$, we define $I_N(y)$ as the first integer $j$ for which

$$\left|\sum_{k=1}^{j} X_k(\omega)k^{-t/2}T^{\mathbf{j}_k}f(y)\right| = \max_{1\leq n\leq N}\left|\sum_{k=1}^{n} X_k(\omega)k^{-t/2}T^{\mathbf{j}_k}f(y)\right|$$

and obtain

$$\lim_{N\to\infty}\left\|\max_{1\leq n\leq N}\left|\sum_{k=1}^{n} X_k(\omega)k^{-(2-q)/(2q)}T^{\mathbf{j}_k}f\right|\right\|_q$$
$$= \lim_{N\to\infty}\|\|\Psi_{I_N,t}(T)f\|\|_q \leq C_{1,2,q}\|f\|_q.$$

This yields the $L_q(\pi)$-integrability of the maximal function. Furthermore, it yields that $\sup_{n\geq 1}|\sum_{k=1}^{n} X_k(\omega)k^{-(2-q)/(2q)}T^{\mathbf{j}_k}f| < \infty$ $\pi$-a.e. for any $f \in L_q(\pi)$. Note that $\tilde{\Omega}^*$ is a subset of the set $\Omega^*$ which was defined in Theorem 5.1. So, for $f \in L_2(\pi)$ the series $\sum_{n=1}^{\infty} X_n(\omega)T^{\mathbf{j}_n}f$ converges $\pi$-a.e.; hence also, by Abel's summation by parts, $\sum_{n=1}^{\infty} X_n(\omega)n^{-(2-q)/(2q)}T^{\mathbf{j}_n}f$ converges $\pi$-a.e. Since $L_2$ is dense in $L_q$, the Banach principle yields that for any $f \in L_q(\pi)$ the series $\sum_{n=1}^{\infty} X_n(\omega)n^{-(2-q)/(2q)}T^{\mathbf{j}_n}f$ converges $\pi$-a.e. $\square$

REMARK. In (15) we could sharpen the estimation in order to obtain a slightly better rate in the normalization $n^{(2-q)/(2q)}$ that appears in Theorems 5.3 and 5.4 below.

Recall that a Dunford–Schwartz operator on $L_1(Y,\pi)$ is a contraction $T$ which is also a contraction of $L_\infty(Y,\pi)$, and therefore is also a contraction of each $L_q(Y,\pi)$, $1 < q < \infty$, by the Riesz–Thorin theorem (for a simple proof for Markov operators, see [21], page 65). Clearly, the operator $T$ defined in the previous theorem is a special kind of Dunford–Schwartz operator.

THEOREM 5.4. *Let $\{X_n\} \subset L_2(\Omega,\mu)$ be a sequence of centered independent random variables, and let $\{p_n\}$ and $\{q_n\}$ be sequences of natural numbers. If the series*

$$\sum_{n=1}^{\infty} \|X_n\|_2^2 \log(n \vee p_n^* \vee q_n^*)(\log n)^2$$

*converges, then there exists a set of full measure $\Omega^* \subset \Omega$, such that for every $\omega \in \Omega^*$, for every commuting Dunford–Schwartz operator $T_1$ and $T_2$ on a space $L_1(Y,\pi)$ and any $f \in L_q(Y,\pi)$, $1 < q \leq 2$, the series*

$$\sum_{n=1}^{\infty} \frac{X_n(\omega)T_1^{p_n}T_2^{q_n}f}{n^{(2-q)/(2q)}} \qquad \textit{converges } \pi\textit{-a.e.}$$

*Furthermore, for every $\omega \in \Omega^*$ there exists a constant $K_\omega < \infty$, independent of $f$, such that*

$$\left\| \sup_{n \geq 1} \left| \sum_{k=1}^{n} \frac{X_k(\omega) T_1^{p_n} T_2^{q_n} f}{n^{(2-q)/(2q)}} \right| \right\|_q \leq K_\omega \|f\|_q.$$

PROOF. In order to apply Stein's interpolation theorem in Theorem 5.3, it was needed that the operators involved there are Dunford–Schwartz. Since the proof of Theorem 5.3 uses Theorem 5.1, the failure of the dilation theorem for more than two commuting contractions restricts the present theorem to the case of only two commuting Dunford–Schwartz operators. □

**6. On the Wiener–Wintner property.** In a series of papers (see [2], the book [3] and the references therein), Assani introduced the concepts of Wiener–Wintner functions and of Wiener–Wintner (dynamical) systems. These families are connected to several deep theorems of Bourgain (e.g., the return times theorem and the double recurrent theorem). We show that the estimates obtained for random trigonometric polynomials allow one to deduce the Wiener–Wintner property for functions (on a dynamical system), which satisfy some appropriate mixing conditions.

Let us recall the notion of Wiener–Wintner functions.

DEFINITION 6.1. Let $(\Omega, \mathcal{F}, \mu, \theta)$ be a dynamical system. For $0 < \alpha < 1$, a function $f$ is a WW function of power type $\alpha$ in $L_p(\Omega, \mu)$, $1 \leq p < \infty$, if there exists a constant $C_f > 0$ such that

$$\left\| \max_{t \in [-\pi, \pi)} \left| \frac{1}{n} \sum_{k=1}^{n} e^{ikt} f \circ \theta^k \right| \right\|_p \leq \frac{C_f}{n^\alpha} \qquad \text{for every } n \geq 1.$$

Now, we have the following corollary of Theorem 3.1.

PROPOSITION 6.1. *Let $(\Omega, \Sigma, \mu, \theta)$ be a dynamical system, and put $\mathcal{F}_k = \sigma\{f, \ldots, f \circ \theta^k\}$. For every $p \geq 1$, there exists $C_p > 0$, such that for any $f \in L_\infty(\Omega, \mu)$,*

$$(16) \quad \left\| \max_{t \in [-\pi, \pi)} \left| \frac{1}{n} \sum_{k=1}^{n} e^{ikt} f \circ \theta^k \right| \right\|_p$$
$$\leq \frac{C_p \sqrt{\log n}}{n} \left( n \|f\|_\infty^2 + \sum_{i=2}^{n} \sum_{k=1}^{i-1} \|f \circ \theta^k \mathbf{E}[f \circ \theta^i | \mathcal{F}_k]\|_\infty \right)^{1/2}.$$

*In particular, if $\{f \circ \theta^n\}$ is a martingale difference sequence, then for every $0 < \alpha < 1/2$ the function $f$ is a WW function of order $\alpha$ in all $L_p(\Omega, \mu)$ spaces.*

Proposition 6.1 requires *uniform* estimation of the correlation coefficients, which may look too restrictive.

It is possible to use Theorem 3.4, in order to obtain:

PROPOSITION 6.2. *Let $(\Omega, \Sigma, \mu, \theta)$ be a dynamical system. For every $p > 2$, there exists $C_p > 0$, such that for every $f \in L_p(\Omega, \mu)$*

$$\left\| \max_{t \in [-\pi, \pi)} \left| \frac{1}{n} \sum_{k=1}^{n} e^{ikt} f \circ \theta^k \right| \right\|_p$$
$$\leq \frac{C_p}{n^{1-1/p}} \left( n \|f\|_p^2 + \sum_{i=2}^{n} \sum_{k=1}^{i-1} \|f \circ \theta^k \mathbf{E}[f \circ \theta^i | \mathcal{F}_k]\|_{p/2} \right)^{1/2}.$$

*In particular, if $\{f \circ \theta^n\}$ generates a martingale difference sequence, then $f$ is a WW function of power type $\alpha = 1 - 1/p - 1/2$ in $L_p(\Omega, \mu)$.*

Clearly, control of the correlation coefficients involved in (16) (or in Proposition 6.2) yields the Wiener–Wintner property. Such a control is possible in many situations. See, for example, the discussion in page 9 of [13]. A typical example is provided by a Markov chain whose transition probability induces a quasi-compact operator (see [19] and the references therein). For results without the quasi-compactness assumption one can refer to [9] or [10].

We now give an application of Proposition 6.1 to $K$-automorphisms.

DEFINITION 6.2 ([36]). Let $(\Omega, \mathcal{F}, \mu)$ be a probability space, and let $\theta$ be an invertible measure-preserving point transformation on $\Omega$. The dynamical system $(\Omega, \mathcal{F}, \mu, \theta)$ is called a $K$-automorphism if there exists a sub-$\sigma$-algebra $\mathcal{C}$, such that

$$\theta^{-1}\mathcal{C} \subset \mathcal{C}; \qquad \bigcap_{n \geq 1} \theta^{-n}\mathcal{C} = \{\varnothing, \Omega\}; \qquad \bigcup_{n \geq 1} \theta^n \mathcal{C} \text{ is dense in } \mathcal{F}.$$

PROPOSITION 6.3. *Let $(\Omega, \mathcal{F}, \mu, \theta)$ be a $K$-automorphism, and let $1 \leq p < \infty$. There exists a set of functions, which is dense in $L_p^0(\mu) = \{f \in L_p(\mu) : \mathbf{E}(f) = 0\}$, such that for every $f$ in this set, there exists a constant $C_{f,p}$ such that*

$$(17) \quad \left\| \max_{t \in [-\pi, \pi]} \left| \sum_{k=1}^{n} e^{int} f \circ \theta^n \right| \right\|_p \leq C_{f,p} \sqrt{n \log(n+1)} \qquad \text{for every } n \geq 1.$$

PROOF. Let $\mathcal{C}$ be the sub-$\sigma$-algebra related to the $K$-automorphism $(\Omega, \mathcal{F}, \mu, \theta)$. By Definition 6.2 the algebra $\bigcup_{n \geq 1} \theta^n \mathcal{C}$ is dense in $\mathcal{F}$, hence (by basic measure theory) the set of functions $c.l.m.\{\mathbf{1}_C - \mu(C) : C \in \bigcup_{n \geq 1} \theta^n \mathcal{C}\}$ is dense in each $L_p^0$ for the norm $\|\cdot\|_p$. By the second requirement in Definition 6.2, as $k$ goes to infinity, the martingale $\mathbf{E}(\mathbf{1}_C|\theta^{-k}\mathcal{C})$ converges a.e. to $\mu(C)$ (see [14], Chapter VII, Theorem 4.3). So by the bounded convergence theorem, $\mathbf{1}_C - \mathbf{E}(\mathbf{1}_C|\theta^{-k}\mathcal{C})$ converges in $L_p$ to $\mathbf{1}_C - \mu(C)$. Hence, the set of functions $\{\mathbf{1}_C - \mathbf{E}(\mathbf{1}_C|\theta^{-k}\mathcal{C}) : k \geq 1, C \in \theta^{-k}\mathcal{C}\}$ is dense in $L_p^0$ for the norm $\|\cdot\|_p$, and it is sufficient to prove that (17) holds for functions of this type.

Let $k \geq 1$, and let $C \in \theta^{-k}\mathcal{C}$. Put $f = \mathbf{1}_C - \mathbf{E}(\mathbf{1}_C|\theta^{-k}\mathcal{C})$. Since $0 \leq \mathbf{E}(\mathbf{1}_C|\theta^{-k}\mathcal{C}) \leq 1$ a.e. we have $|f| \leq 1$ a.e. Since $\theta^{-j}\mathcal{C} \subset \theta^{-k}\mathcal{C}$ for every $j \geq k$, we have

$$\mathbf{E}(f|\theta^{-j}\mathcal{C}) = \mathbf{E}(\mathbf{1}_C|\theta^{-j}\mathcal{C}) - \mathbf{E}(\mathbf{E}(\mathbf{1}_C|\theta^{-k}\mathcal{C})|\theta^{-j}\mathcal{C}) = 0 \tag{1}$$

for every $j \geq k$.

CLAIM. *Let $\Sigma$ be a sub-$\sigma$-algebra of $\mathcal{F}$, and let $\eta$ be an invertible measure-preserving transformation. Then for any integrable random variable $Z$, we have $\mathbf{E}(Z \circ \eta | \Sigma) = [\mathbf{E}(Z|\eta\Sigma)] \circ \eta$ a.e.*

PROOF. For any $B \in \Sigma$ we have

$$\int_B \mathbf{E}(Z \circ \eta | \Sigma)\, d\mu = \int \mathbf{1}_B \cdot Z \circ \eta\, d\mu = \int (\mathbf{1}_{\eta B} \cdot Z) \circ \eta\, d\mu$$

$$= \int \mathbf{1}_{\eta B} \cdot Z\, d\mu = \int_{\eta B} \mathbf{E}(Z|\eta\Sigma)\, d\mu$$

$$= \int [\mathbf{1}_{\eta B} \cdot \mathbf{E}(Z|\eta\Sigma)] \circ \eta\, d\mu = \int_B [\mathbf{E}(Z|\eta\Sigma)] \circ \eta\, d\mu.$$

Since $[\mathbf{E}(Z|\eta\Sigma)] \circ \eta$ is $\Sigma$-measurable, the result follows from the uniqueness of the conditional expectation. □

For every $1 \leq i \leq n$, put $X_i = f \circ \theta^{n+1-i}$, and for $i > n$ put $X_i \equiv 0$. Since $|f| \leq 1$, also $|X_i| \leq 1$. As usual, for any $l \geq 1$ put $\mathcal{F}_l = \sigma(X_1, \ldots, X_l)$. Since $f$ is $\mathcal{C}$-measurable, $X_i$ is $\theta^{-(n+1-i)}\mathcal{C}$-measurable, so $\mathcal{F}_l \subset \theta^{-(n+1-l)}\mathcal{C}$. Let $1 \leq j \leq l$ with $l - j \geq k$. Using Claim and equality (1) we have

$$\mathbf{E}(X_l|\mathcal{F}_j) = \mathbf{E}(\mathbf{E}(f \circ \theta^{n+1-l}|\theta^{-(n+1-j)}\mathcal{C})|\mathcal{F}_j)$$

$$= \mathbf{E}([\mathbf{E}(f|\theta^{-(l-j)}\mathcal{C})] \circ \theta^{n+1-l}|\mathcal{F}_j) = 0.$$

With our previous notation, we obtain

$$R_{0,n} = \sum_{i=1}^n \|X_i\|_\infty^2 + \sum_{i=1}^n \sum_{j=1}^i \|X_j \mathbf{E}(X_i|\mathcal{F}_j)\|_\infty$$

$$= \sum_{i=1}^n \|X_i\|_\infty^2 + \sum_{j=1}^n \sum_{i=j}^n \|X_j \mathbf{E}(X_i|\mathcal{F}_j)\|_\infty \leq n + nk.$$

By Proposition 6.1, there exists some universal constant $C_p$, such that

$$\left\|\max_{t\in[-\pi,\pi]}\left|\sum_{j=1}^n e^{ijt}X_j\right|\right\|_p \leq C_p\sqrt{kn\log(n+1)}.$$

Then it follows easily that

$$\left\|\max_{t\in[-\pi,\pi]}\left|\sum_{j=1}^n e^{ijt}f\circ\theta^j\right|\right\|_p \leq \sqrt{k}C_p\sqrt{n\log(n+1)},$$

where $C_{f,p} = \sqrt{k}C_p$. $\square$

REMARKS. 1. The estimate (17) considerably improves the estimate obtained in the proof of Theorem 4 in [2]. Similarly, our results can be applied to improve Theorems 6 and 7 of [2].

2. The dense set of functions appearing in Proposition 6.3 is the same as the one considered in [2].

3. The constant $C_{f,p}$ may be chosen uniformly for $\{f\circ\theta^l : l\in \mathbf{Z}\}$. Hence the estimate (17) can be obtained along blocks, and we may apply Proposition 2.6 (see also the remarks after it) in order to deduce uniform convergence of the one-sided rotated Hilbert transform, even with rate.

4. As noticed by Assani, Proposition 6.3 gives an example of a Wiener–Wintner dynamical system of all powers $0 < \alpha < 1/2$.

**Acknowledgments.** The manuscript was completed during the first author's postdoctoral fellowship at the E. Schröedinger Institute, Vienna. The first author is very grateful to Paul Furhmann for his advice and encouragement. Both authors are very grateful to Michael Lin for mathematical discussions and his constant enthusiasm.

Department of Mathematics
Ben Gurion University
P.O.B. 653
841105 Beer Sheva
Israel
e-mail: guycohen@ee.bgu.ac.il

Department of Mathematics
University of New Caledonia
Equipe ERIM
B.P. 4477
F-98847 Noumea Cedex
France
e-mail: cuny@univ-nc.nc